\documentclass[pdflatex,sn-mathphys-ay]{sn-jnl}

\usepackage{graphicx}%
\usepackage{multirow}%
\usepackage{amsmath,amssymb,amsfonts}%
\usepackage{amsthm}%
\usepackage{mathrsfs}%
\usepackage[title]{appendix}%
\usepackage{xcolor}%
\usepackage{textcomp}%
\usepackage{manyfoot}%
\usepackage{booktabs}%
\usepackage{algorithm}%
\usepackage{algorithmicx}%
\usepackage{algpseudocode}%
\usepackage{listings}%

\theoremstyle{thmstyleone}%
\newtheorem{theorem}{Theorem}[section]

\newtheorem{lemma}[theorem]{Lemma}%

\newtheorem{corollary}[theorem]{Corollary}%

\theoremstyle{thmstyletwo}%
\newtheorem{example}{Example}%
\newtheorem{remark}{Remark}%

\theoremstyle{thmstylethree}%
\newtheorem{definition}{Definition}%

\usepackage{hyperref}
\usepackage{graphicx}
\usepackage{xcolor}
\usepackage{subcaption}
\usepackage{bm}

\usepackage[english]{babel}
\addto\captionsenglish{\renewcommand{\appendixname}{Appendix}}
\usepackage{soul}
\usepackage{lineno}
\usepackage{comment}

\numberwithin{equation}{section}

\DeclareMathOperator*{\diag}{diag}
\DeclareMathOperator*{\tr}{{\rm tr}}
\DeclareMathOperator*{\E}{{\rm E}}

\DeclareMathOperator*{\argmin}{arg\,min}

\newcommand{\IR}{\mathbb{IR}}
\newcommand{\R}{\mathbb{R}}
\newcommand{\Var}{{\rm Var}}
\newcommand{\Cov}{{\rm Cov}}

\newcommand{\U}{U}
\newcommand{\Uc}{\tilde{U}}

\graphicspath{{Fig/}}

\begin{document}

\title[Location and association measures for interval-valued data based on Mallows' distance]{Location and association measures for interval-valued data based on Mallows' distance}

\author*[1,2]{\fnm{M. Ros\'ario} \sur{Oliveira}}\email{rosario.oliveira@tecnico.ulisboa.pt}

\author[2]{\fnm{Diogo} \sur{Pinheiro}}\email{diogo.pinheiro.99@tecnico.ulisboa.pt}
\equalcont{These authors contributed equally to this work.}

\author[2,3]{\fnm{Lina} \sur{Oliveira}}\email{lina.oliveira@tecnico.ulisboa.pt}
\equalcont{These authors contributed equally to this work.}

\affil*[1]{\orgdiv{Department of Mathematics}, \orgname{Instituto Superior T\'ecnico},  \city{Lisbon},  \country{Portugal}}

\affil[2]{\orgdiv{CEMAT},  \orgname{Instituto Superior T\'ecnico}, \city{Lisbon}, \country{Portugal}}

\affil[3]{\orgdiv{CAMGSD}, \orgname{Instituto Superior T\'ecnico},  \city{Lisbon}, \country{Portugal}}


\abstract{
The growing demand to analyse large and complex datasets has spurred the development of Symbolic Data Analysis as a promising approach to address contemporary data challenges. Amongst these, interval-valued data introduces new theoretical and methodological questions that remain open.

In this paper, we generalise measures of location and association for interval-valued random variables using Mallows’ distance. Departing from restrictive assumptions such as uniform distributions over microdata, our proposal extends the barycentre approach to any absolutely continuous distribution with finite second moment. A key contribution is the derivation of explicit formulas for Mallows’ distance in \textit{p}-dimensional interval spaces. These formulas decompose into components for centres, ranges, and a novel cross-term that captures their interaction. This decomposition leads to a new theoretical symbolic covariance matrix that explicitly accounts for the dependence between centres and ranges — a relation often obscured in current definitions of symbolic covariance.

Theoretical developments are supported by empirical studies on diverse real-world datasets, each reflecting different degrees of information about the underlying microdata. These applications highlight both the flexibility of the proposed methodology and the interpretability of its results.}


\maketitle

\section{Introduction}\label{Sec:Intro}
The explosion of data in recent decades has motivated the emergence of new data types and the demand for more complex statistical techniques to address them. 
Symbolic data analysis (SDA) is a field of statistics that studies data with internal variation, of which histograms and intervals are two key examples. 
Symbolic objects appear mainly from the agglomeration of individual real-valued or categorical observations. SDA relies on statistical methods to learn patterns from individuals (microdata) based on aggregate observations (macrodata). 
Common situations for using symbolic data include a large sample size, privacy concerns, research interests, and the symbolic nature of the data being collected. 
For a thorough review of symbolic data types and their analysis, see \cite{BillardDiday2006,billard.diday:2020,Brito:2014}.

The works \cite{Bertrand.Goupil:2000,BillardDiday2006,Billard2008} introduced measures of location, dispersion, and association between interval-valued random variables, formalised as a function of the observed macrodata and implicit assumptions about the microdata. 
\cite{Bertrand.Goupil:2000} proposed that the sample mean and sample variance of a set of interval-valued observations should simply be their centres' sample mean and sample variance, respectively. 
In \cite{Irpino.verde:2015} this approach was called ``\textit{SDA two-level paradigm}''. The authors suggested an alternative where the location measure was the Fréchet mean, or barycentre, of the set of symbolic observations. They considered the space of real bounded intervals and the $L_2$ Wasserstein distance, also known as the Mallows' distance (denomination used in the rest of this paper), based on the assumption that the microdata spread in each observed interval according to a uniform distribution. Under this approach, the location measure is an interval, by contrast with the previous definition of this measure as a real number, and the variance is a non-negative real number, as usual. In this paper, we generalise the barycentre approach to the population framework, admitting any possible absolutely continuous distribution with finite second moment for the microdata. 

The sample covariance and sample correlation matrices were also addressed in the context of symbolic principal component analysis in \cite{rademacher2012symcovpca,Oliveira.et.al:2017,wang2012cipca}. Specifically, in \cite{Oliveira.et.al:2017} the authors established relationships between several proposed methods of symbolic principal component analysis and available definitions of sample symbolic variance and covariance. Later, in \cite{Serrao.et.al:2023} the principal components were derived as the linear combinations of the original interval-valued random variables which maximised the symbolic variance.

Other areas of Statistics have also been addressed by SDA, like
clustering 
\citep{billard.diday:2020,Carvalho.Lechevallier:2009,Sato.Ilic:2011},
discriminant analysis 
\citep{Dias.Brito:2021,Queiroz.et.al:2018,DSilva.Brito:2015},
regression analysis 
\citep{LNeto.et.al:2011,Dias.Brito:2017,irpino.Verde:2015LReg,Whitaker.et.al:2021}, 
time series 
\citep{Lin.GonzalesRivera:2016,Maia.et.al:2008,Teles.Brito:2015},
Bayesian hierarchical modelling 
\citep{Lin.et.al:2022}, and
network sciences 
\citep{Alves.et.al:2022,Ponti.Irpino.et.al:2022}), amongst others.

Parametric approaches for interval-valued variables have also been considered. In \cite{LeRademacher.Billard:2011}, the authors derived maximum likelihood estimators for the mean and the variance of interval-valued and histogram-valued variables, assuming  uniform or symmetric triangular distributions for the microdata. In the follow-up paper \cite{Samadi.Billard.et.all:2023} revised the initial work and derived the maximum likelihood estimators for all important covariance statistics.
In \cite{LNeto.et.al:2011}, interval-valued variables were formulated as bivariate random vectors to introduce a symbolic regression model based on the theory of generalised linear models. The contributions in \cite{Brito2012,DSilva.Brito:2015,Duarte.Silva.Brito:2018} followed a different approach. The centres and logarithms of the ranges were collected in a random vector with a multivariate normal or skew-normal distribution, which was used to derive methods for the analysis of variance \citep{Brito2012}, discriminant analysis \citep{DSilva.Brito:2015}, and outlier detection \citep{Duarte.Silva.Brito:2018} of interval-valued variables.
More recently, a line of research was developed using likelihood-based methods that fitted models for the microdata when only the macrodata were observed \citep{Rahman.et.al.2022,Zhang.Sisson:2020}.

In this paper, we consider the interval-valued data model that establishes the connection between macrodata and microdata, as introduced in \cite{Oliveira.et.al:2021}. In Section~\ref{Sec:Background}, we present the model along with the necessary background. In Section~\ref{Sec:Mallows}, we derive general expressions for the Mallows’ distance between two $p$-dimensional intervals, expressed in terms of their centers, ranges, and the first two moments of the corresponding microdata distributions. We also demonstrate that this Mallows' distance can be interpreted as a special case of a weighted Euclidean distance (see Section~\ref{Sec:WasserstIRp}). 
These results provide the theoretical foundation for defining location, scale, and association measures of interval-valued data using the barycentre approach, discussed in Section~\ref{Sec4:moments}.
We illustrate our proposal in Section~\ref{Sec:Examples} through three examples, reflecting different levels of available microdata information. We also discuss the selection of distributions and parameter values to model real datasets.
In Section~\ref{Sec:concl}, we present the main conclusions. The proofs for the results in Section~\ref{Sec:Mallows} and Section~\ref{Sec4:moments} can be found in Appendix~\ref{AppendixA} and Appendix~\ref{AppendixB}, respectively.


\section{Preliminaries}\label{Sec:Background}

We motivate interval-valued data with an example. 
During 2024, a weather station located in Quinta da França, Portugal, monitored the local temperature every five minutes. Each day, the station recorded 288 values. To compare the daily temperatures over different days, it is common to use one or more summary statistics, such as the mean and the standard deviation. Although this is practical, it leads to a loss of information on the distribution of the data. 
Our symbolic approach proposes to aggregate the daily temperatures into a single real-valued interval onto which the distribution of the original data is attached. This symbolic construction, known as a \textit{symbolic interval}, consists of two hierarchically linked components. 
The real-valued interval as a set of points between two real numbers is the \textit{macrodata} and the individual points are the \textit{microdata}.

\begin{figure}[h] 
\hspace{10pt}
\begin{subfigure}[t]{0.3\textwidth}
\centering
\subcaptionbox{6th March: [-1.3,13.2]}{
\includegraphics[alt = {Temperature Interval of March 6th. Two humps located at the endpoints}, width=1.1\textwidth]{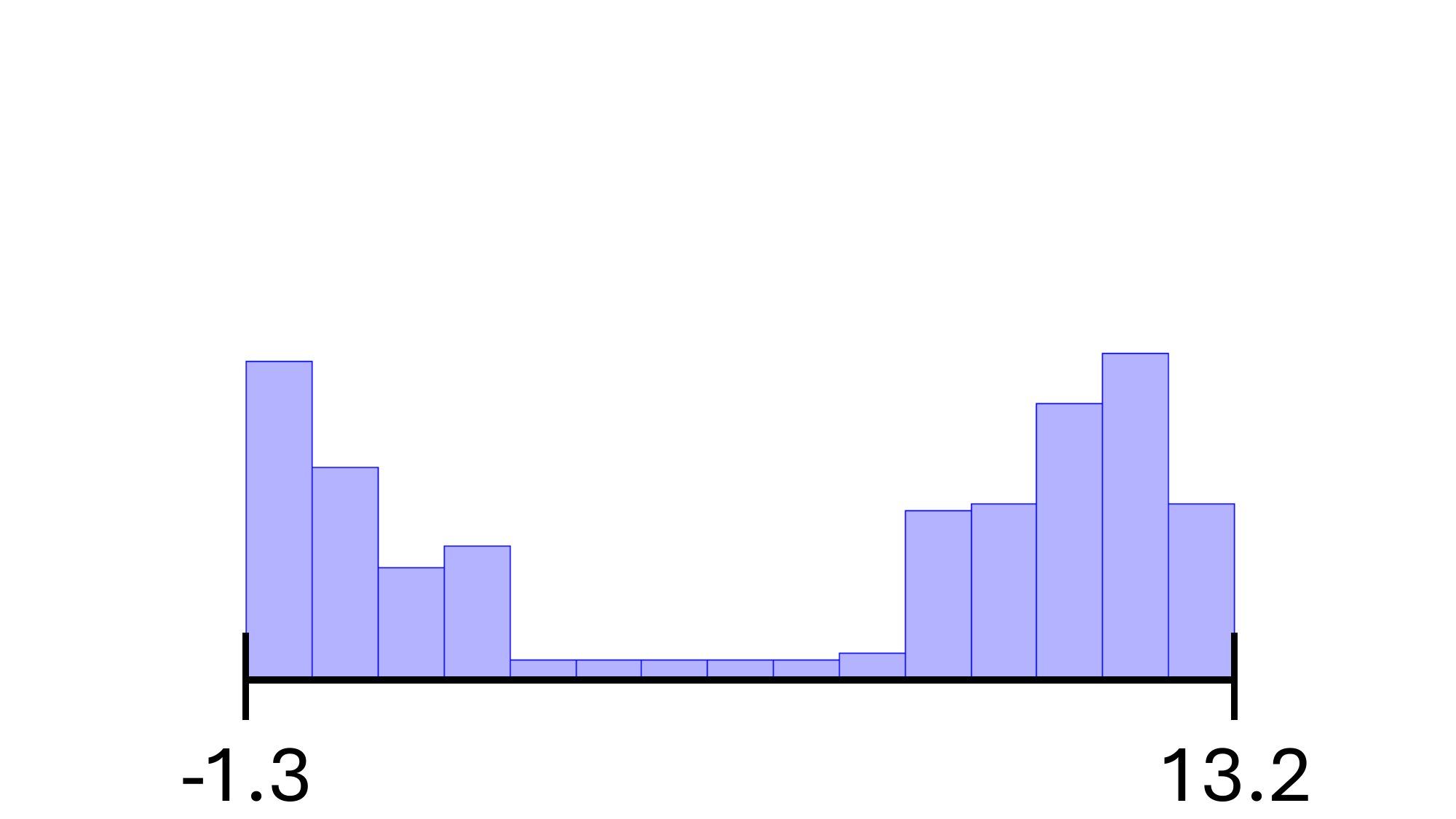}}
\label{Fig:QuintaDay6Mar}
\end{subfigure}
\begin{subfigure}[t]{0.3\textwidth}
\centering
\subcaptionbox{14th July: [9.1,28.6]}{\includegraphics[alt = {Temperature Interval of July 14th. Higher temperatures overall with peaks in endpoints}, width=1.1\textwidth]{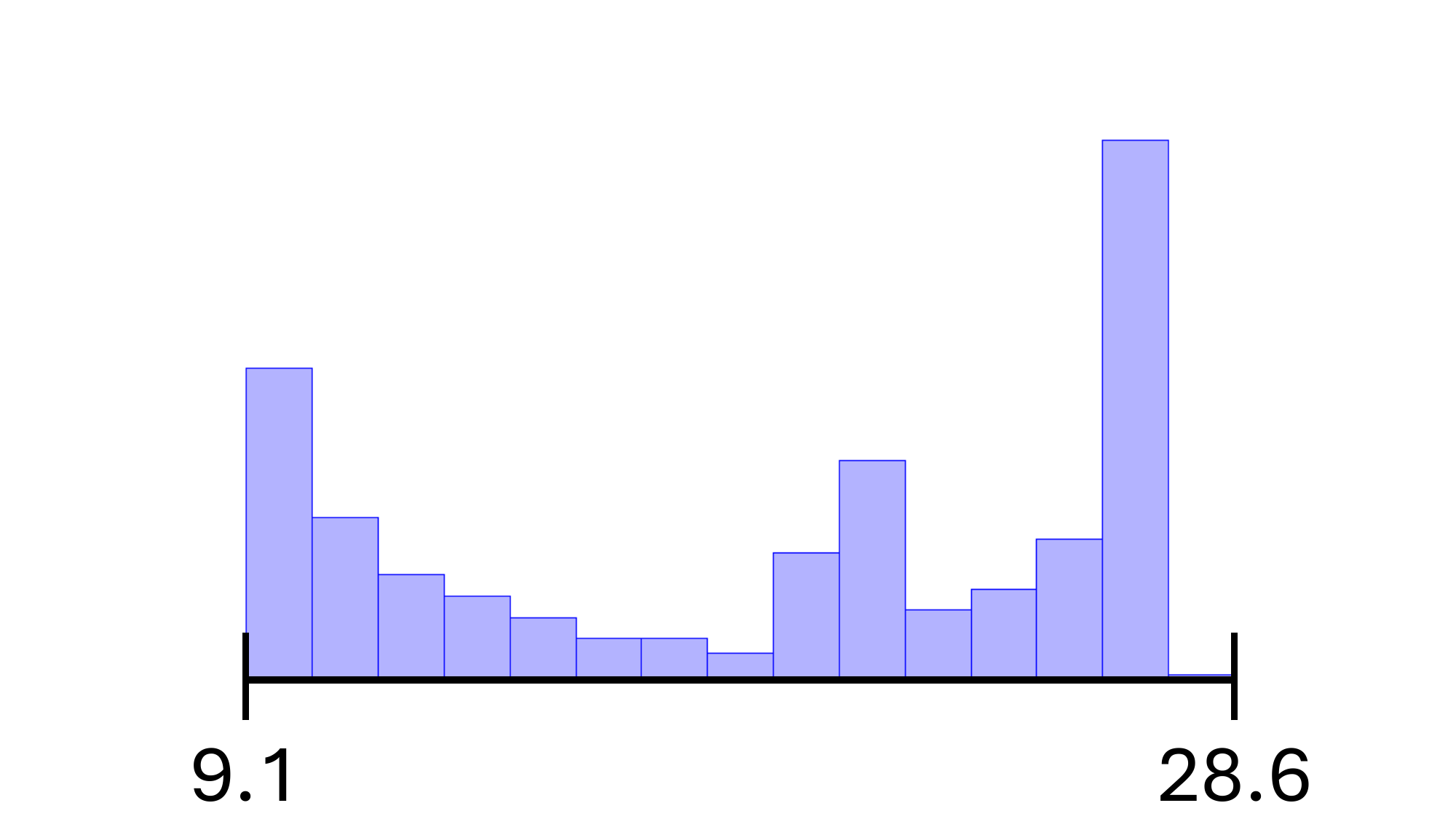}}
\label{Fig:QuintaDay14Jul}
\end{subfigure}
\begin{subfigure}[t]{0.3\textwidth}
\centering
\subcaptionbox{6th December: [4.6,19.2]}{\includegraphics[alt = {Temperature Interval of December 6th reveals a skewed distribution to the lower endpoint}, width=1.1\textwidth]{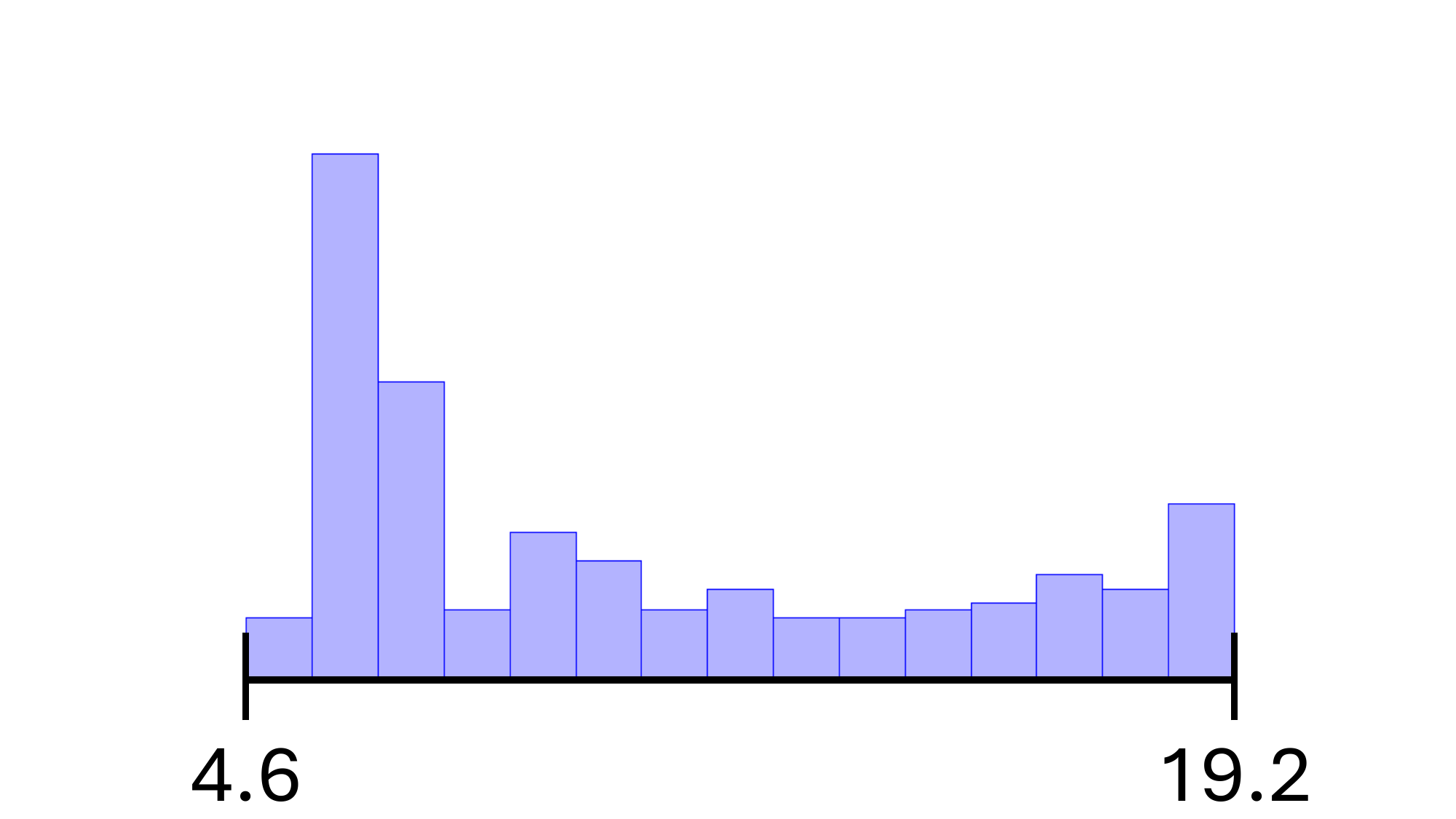}}
\label{Fig:QuintaDay6Dec}
\end{subfigure}

\caption{Examples of symbolic intervals of daily temperatures in Quinta da França, Portugal}
\label{fig:QuintaDayInt}
\end{figure}

Figure \ref{fig:QuintaDayInt} illustrates symbolic intervals created using daily temperatures recorded on three different days. Observing the macrodata intervals, we can see that the minimum and maximum temperatures vary substantially depending on the season of the year. More information can be extracted from the microdata. We use histograms to depict the frequency distribution of the data. The tendency to have a larger mass in the endpoints of the intervals is a consequence of the day and night cycle. Skeweness in the distribution is also a term of comparison and indication of the season.

We can represent the symbolic interval by the tuple $x = ([a,b], F)$, where $[a,b] \subset \R$ is the real-valued interval of macrodata and $F$ is an absolutely continuous distribution function with support $[a,b]$. 
Note that this function can be estimated when we have access to the microdata. In most cases, however, we only observe the macrodata and the distribution function of the microdata needs to be assumed. The continuous uniform distribution is the most common assumption in the literature.

\begin{definition}\label{IRp}
Let $\IR=\{[a,b]\colon a,b\in\R,\, a\leq b\}$ be the set of all real closed bounded intervals. For a positive integer $p$, let $\IR^p$ be the cartesian product of $p$ copies of $\IR$, that is,
  \begin{equation*}
        \IR^p = \{([a_1, b_1], \ldots, [a_p, b_p])^T:\,  a_i, b_i \in \mathbb{R}, \,a_i \leq b_i,\, i=1,\dots, p\},
\end{equation*}
the set of closed bounded $p$-dimensional intervals or hyperrectangles.
\end{definition}

A real-valued interval $[a,b] \in \IR$ is uniquely defined by its \emph{centre} ${c= (a+b)/2} \in \R$ and \emph{range} $r=b-a \in \R^{+}_0$. 
Hence, there exists a correspondence between $\IR$ and $\R \times \R_0^+$ through the bijective mapping which sends the interval $[a,b]$ to $(c,r)^T\in \R \times \R_0^+$.
Note that this correspondence extends to higher dimensions. Given a real-valued hyperrectangle 
$([a_1, b_1], \ldots, [a_p, b_p])^T \in \IR^p$, $p\in\mathbb{N}$, let
$ (\bm{c}^T, \bm{r}^T)^T \in \R^p\times (\R^+_0)^p$ be the vector of centres and ranges, that is,
\begin{equation*}
\bm{c}=\left(\frac{a_1 + b_1}{2}, \ldots, \frac{a_p + b_p}{2} \right)^T \quad {\rm and } \quad \bm{r} = (b_1 - a_1, \ldots, b_p - a_p)^T.
\end{equation*}

Using a slight abuse of notation, the bijection between $\IR^p$ and $\R^p\times (\R^+_0)^p$ allows us to simplify and write
\begin{equation*}
([a_1, b_1], \ldots, [a_p, b_p])^T=(\bm{c}^T, \bm{r}^T)^T.
\end{equation*}

As suggested in Figure \ref{fig:QuintaDayInt}, the macrodata can vary substantially between individuals. This adds to the difficulty of studying the microdata and their distribution functions with different supports. We circumvent this issue by normalising the microdata to the interval $[-1,1]$ through a linear transformation that removes the contributions of the centre and the range. In more detail, consider the symbolic interval $x = ([a,b], F)$ with centre $c = (a+b)/2$ and range $r = b-a$, and an absolutely continuous real-valued random variable $V$ with support $[a,b] = [c-r/2, c+r/2]$ and distribution function $F$. Furthermore, let $U$ be an absolutely continuous real-valued random variable with support $[-1,1]$, such that the transformation $V = c + r~U/2$ holds. Here, $U$ plays the role of a latent random variable describing the normalised microdata. It is important to mention that the distribution of $U$ becomes clear if $F$ belongs to the location-scale family of distributions.

The previous transformation proposes that the symbolic interval $x$ is completely identified by the centre $c$, the range $r$, and the distribution function of the latent random variable ${U}$ with support $[-1,1]$, say $F_U$. Therefore, we can introduce a more practical notation and state that $x = \left(c,r,F_{U}\right)$, where $c$ and $r$ refer to the macrodata and $F_{U}$ characterises the microdata. This notation is also useful to define interval-valued random variables.

\begin{definition}
Let $\bm{A} = (A_1, \ldots, A_p)^T$ and $\bm{B} = (B_1, \ldots, B_p)^T$ be real-valued random vectors where ${\rm P}(A_i \leq B_i) = 1$, $i = 1,\ldots,p \in \mathbb{N}$. In addition, let $\{U_i, i=1,\ldots,p\}$ be a set of independent and absolutely continuous real-valued random variables with support $[-1,1]$ and distribution functions $F_{U_i}$. We define $\bm{X}=(X_1, \ldots, X_p)^T$ as an interval-valued random vector characterized by $([A_1,B_1],\ldots, [A_p,B_p], F_{U_1}, \ldots, F_{U_p})$.

Alternatively, we write $\bm{X}=(\bm{C},\bm{R}, F_{\bm{U}})$, where $\bm{C}=(C_1, \ldots, C_p)^T$, $\bm{R}=(R_1, \ldots, R_p)^T$ are the real-valued random vectors of centres and ranges, with $\bm{C}=(\bm{B}+\bm{A})/2$, $\bm{R}=\bm{B}-\bm{A}$, and $\bm{U}=(U_1,\ldots,U_p)^T$ is a real-valued random vector with distribution function $F_{\bm U}$.
\end{definition}

In the previous definition, the conventional random vector, $\bm{X}=\bm{C}$, is obtained as a particular case by setting ${\rm P}(R_i=0)=1$ and ${\rm P}(U_i=0)=1$, $i = 1,\ldots,p$. 

It is important to note that interval-valued random variables are typically defined in the literature as the random variables of the macrodata, i.e., $([A_1,B_1],\ldots, [A_p,B_p])^T$. Our proposal expands on this definition by explicitly incorporating the distribution of the microdata, $\bm{X} = ([A_1,B_1],\ldots, [A_p,B_p], F_{\bm{U}})$. This has the advantage of clarifying some gaps related to interval-valued data.

In SDA, macrodata can be seen as the manifest variable that gives information about the behaviour of the microdata that may not be observed. 
The statistical model that establishes a natural link between macrodata and microdata was proposed by 
\cite{Oliveira.et.al:2021}. This model was suggested to unify and add interpretability to definitions of sample interval-valued covariance matrices available in the literature.
The population counterparts of location, scale, and association were also proposed. The model, which has proved its relevance in other areas of SDA \citep{Serrao.et.al:2023,Pinheiro.et.al:2025}, is presented in the next definition.

\begin{definition}\label{Def:Model_Micro_Macro}
    Let $\bm{X}=(\bm{C},\bm{R}, F_{\bm{U}})$ be an interval-valued random vector, where $\bm{C} = (C_1, \ldots, C_p)^T$ and $\bm{R} = (R_1, \ldots, R_p)^T$ are the corresponding random vectors of centres and ranges. Additionally, $\bm{U}=(U_1,\ldots,U_p)^T$ is a random vector of independent and absolutely continuous random variables with support $[-1,1]$.

    The real-valued random vector $\bm{V}=(V_{1},\ldots,V_{p})^T$, with support  $\R^p$, describing the microdata within the macrodata of the interval-valued random vector $\bm{X}$ is defined by
\begin{equation}
\label{Eq:Aij}
	V_{i} = \begin{cases}
	    C_i+U_{i}~\dfrac{ R_i}{ 2},\quad &{\rm if}\quad {\rm P}(R_i=0)=0\\
        C_i,  \quad &{\rm if}\quad {\rm P}(R_i=0)=1\;\wedge\;{\rm P}(U_i=0)=1
	\end{cases}.
\end{equation}
\end{definition}

Observe that, under this model, a realisation of $\bm{V}$ is a point in the random hyperrectangle related to the macrodata of the interval-valued random vector $\bm{X}$, characterised by its centre $\bm{C}$ and range $\bm{R}$. According to the model, the microdata values for a specific real hyperrectangle $[\bm{c}-\bm{r}/2\, ,\, \bm{c}+\bm{r}/2]$, are described by the random vector $\Tilde{\bm{V}} = (\Tilde{V}_1,\ldots,\Tilde{V}_p)^T$, where $\Tilde{V}_i = c_i + \Tilde{U}_ir_i/2$ is the random variable $V_i$ conditioned on realisations of $C_i$ and $R_i$. Likewise, ${\Tilde{U}_i = U_i|(C_i = c_i, R_i = r_i)}$ is the corresponding conditioned latent random variable. 

The assumptions about $U_i$ are typically based on domain knowledge and, at most, are supported by goodness-of-fit measures for methods relying on this formulation. To simplify the proposed model, a reasonable assumption can be made: $\bm{U}$ is independent of the random vector $(\bm{C}^T,\bm{R}^T)^T$. This assumption plays a key role in defining the covariance matrix of an interval-valued random vector and is going to be assumed to hold throughout the remainder of the paper. Under this assumption, ${\Tilde{U}_i = U_i|(C_i = c_i, R_i = r_i)}$ and $U_i$ are identically distributed random variables. Therefore, we write $\bm{X} = (\bm{C}, \bm{R}, F_{\bm{U}})$, with the realisations denoted by $\bm{x} = (\bm{c}, \bm{r}, F_{\bm{U}})$. 
Moreover, $\bm{X}_1,\ldots,\bm{X}_n$ is said to form a random sample of $\bm{X}$ if, for any $i \neq j$, $\bm{X}_i$ and $\bm{X}_j$ are independent and identically distributed interval-valued random vectors; that is, $(\bm{C}_i^T, \bm{R}_i^T, \bm{U}_i^T)^T$ and $(\bm{C}_j^T, \bm{R}_j^T, \bm{U}_j^T)^T$ are independent and identically distributed random vectors. Consequently, each $\bm{X}_j$ is characterized by $(\bm{C}_j, \bm{R}_j, F_{\bm{U}})$ for $j = 1,\ldots,n$. A realization of the random sample is denoted by $\bm{x}_1,\ldots,\bm{x}_n$, or, equivalently, by $(\bm{c}_1,\bm{r}_1,F_{\bm{U}}),\ldots,(\bm{c}_n,\bm{r}_n,F_{\bm{U}})$.

\begin{figure}[h] 
\hspace{68.5pt}
\begin{subfigure}[t]{0.3\textwidth}
\centering
\subcaptionbox{\makebox[110pt][l]{$\text{InvTriang}(-1,1,0)$, $\Var(U_i)=\frac{1}{2}$}}{
\includegraphics[alt = {Probability density function of the symmetric inverted triangular distribution with support in the interval from -1 to 1}, width=0.9\textwidth]{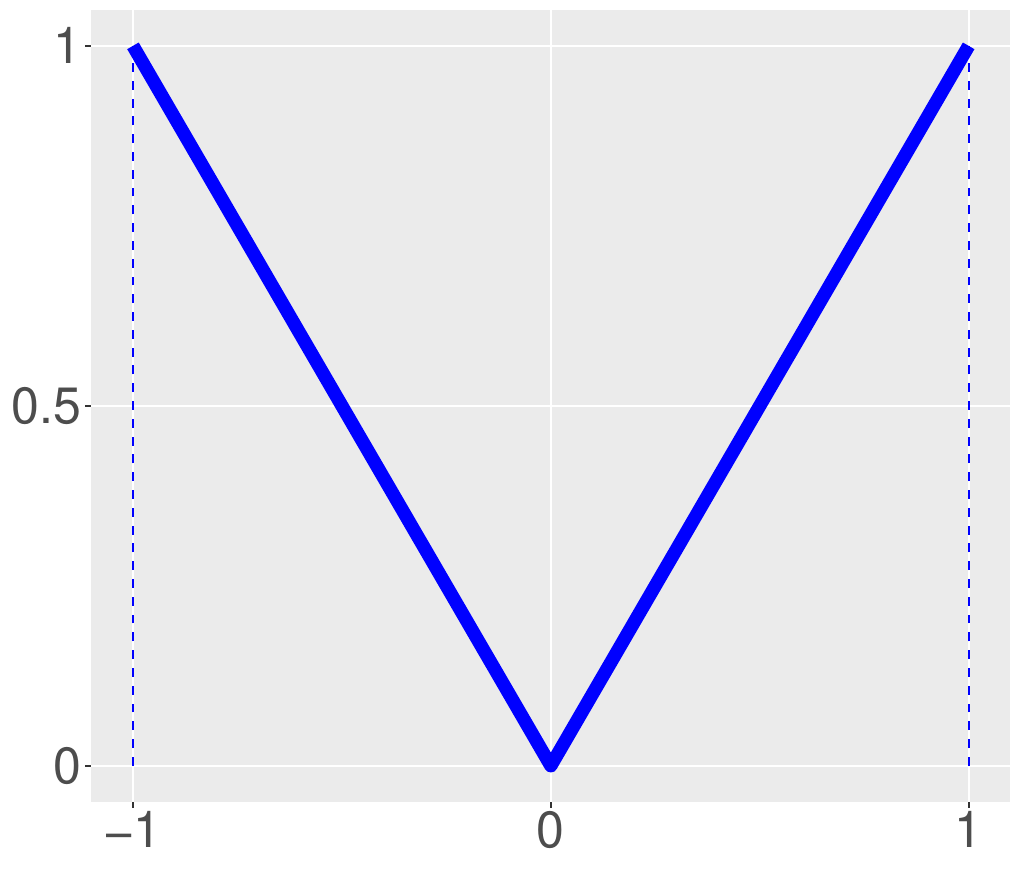}}
\label{Fig:InvT}
\end{subfigure}
\hspace{41.5pt}
\begin{subfigure}[t]{0.3\textwidth}
\centering
\subcaptionbox{\makebox[84pt][l]{$\text{Unif}(-1,1)$, $\Var(U_i)=\frac{1}{3}$.}}{\includegraphics[alt = {Probability density function of the uniform distribution with support in the interval from -1 to 1},width=0.9\textwidth]{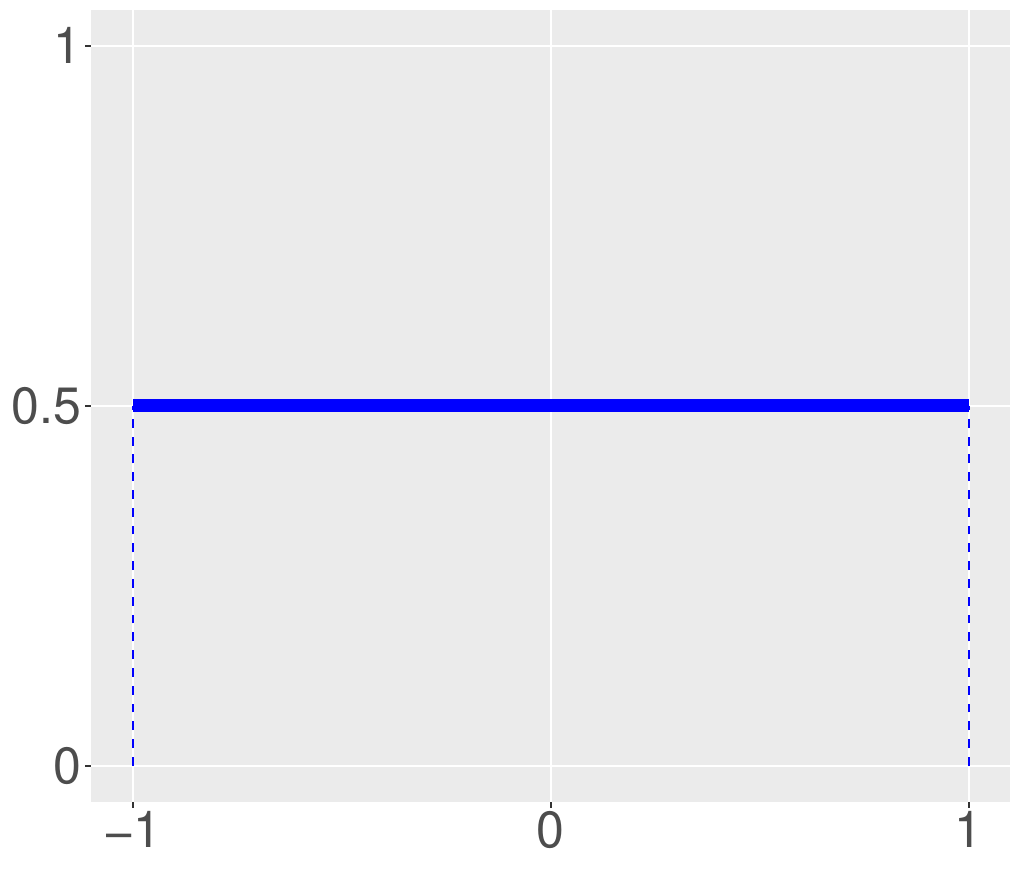}}
\label{Fig:CUnif}
\end{subfigure}

\medskip

\hspace{70pt}
\begin{subfigure}[t]{0.3\textwidth}
\centering
\subcaptionbox{\makebox[100pt][l]{$\text{Triang}(-1,1,0)$, $\Var(U_i)=\frac{1}{6}$.}}
{\includegraphics[alt = {Probability density function of the symmetric triangular distribution with support in the interval from -1 to 1}, width=0.9\textwidth]{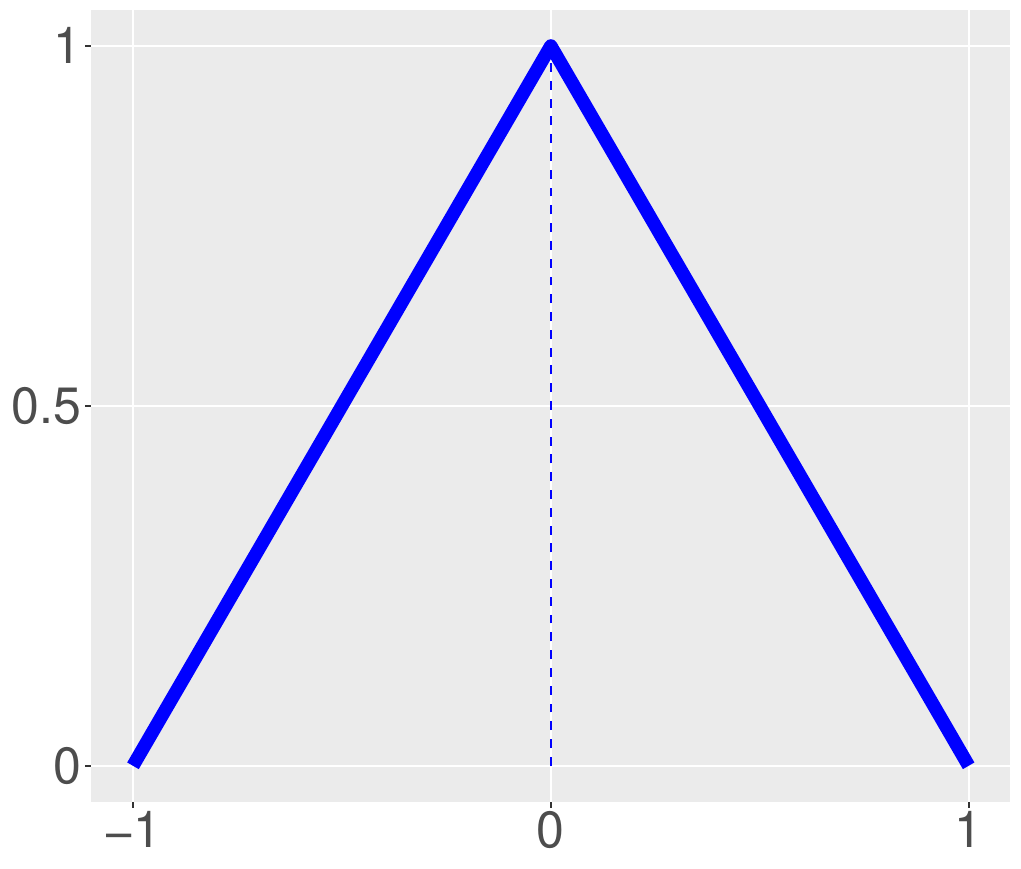}\label{Fig:Tri}}
\end{subfigure}
\hspace{40pt}
\begin{subfigure}[t]{0.3\textwidth}
\centering
\subcaptionbox{\makebox[100pt][l]{$\mathcal{N}\big(0,\frac{1}{9}\big)\big|[-1,1]$, $\Var(U_i)\simeq \frac{1}{9}$.}}
{\includegraphics[alt = {Probability density function of the normal distribution truncated in the interval from -1 to 1, with mean 0 and variance one ninth}, width=0.9\textwidth]{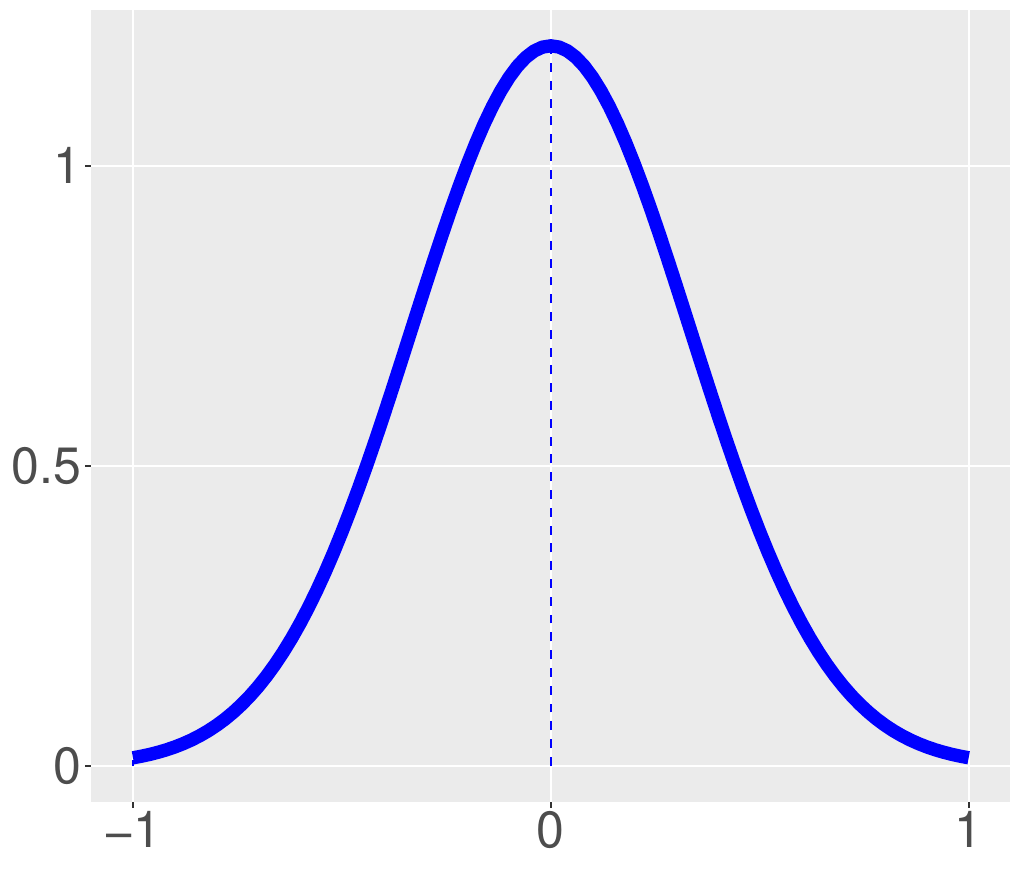}}
\label{Fig:Norm}
\end{subfigure}

\caption{Examples of density functions of continuous symmetric 
distributions of $U_i$ (see Definition~\ref{Def:Model_Micro_Macro})}
\label{fig:distributions}
\end{figure}

In the literature, the most common distribution assumption about the microdata is that they follow a continuous uniform distribution.
However, in \cite{cazes1997cpca} and \cite{Oliveira.et.al:2021} other symmetric alternatives were discussed. Figure~\ref{fig:distributions} illustrates some of those distributions, organised according to the variance of $U_i$, from highest to lowest. In this order, we present the symmetric inverted triangular distribution, InvTriang$(-1,1,0)$, with variance equal to $1/2$, the continuous uniform distribution, Unif$(-1,1)$, with variance equal to $1/3$, the symmetric triangular distribution, Triang$(-1,1,0)$, with variance equal to $1/6$, and the truncated normal distribution, $\mathcal{N}(0,1/9)|[-1,1]$, with variance equal to $1/9 - 2 \phi(3)/(6\Phi(3)-3) \simeq 1/9$, where $\phi(\cdot)$ and $\Phi(\cdot)$ are the probability density function and distribution function of a standard normal distribution, respectively.


\section{Mallows' distance}
\label{Sec:Mallows}

The Mallows' distance has been widely used in SDA to compare intervals and plays an important role in this area (see, for example,  \cite{Brito:2014}). One of the reasons for this is that the Mallows' distance may be seen as a generalisation of the Euclidean distance and shares many of its properties. 
We devote the next two sections to extending and discussing relevant results related with this distance and its role in the definition of location and association for interval-valued objects.

\subsection{Univariate Mallows' distance}\label{Sec:WasserstIR}

We begin by defining the Mallows' distance between two symbolic intervals.

\begin{definition}\label{def:Mallows_unidim}
    Let $x_1=([a_1,b_1],F_1)$ and $x_2=([a_2,b_2],F_2)$, where $F_i$ and $F_i^{-1}$ are, respectively, the distribution function (assumed to have finite second moment) and the quantile function related to the microdata in $[a_{i},b_{i}]$, with centre $c_i=(a_i+b_i)/2$ and range $r_i=b_i-a_i$, $i = 1, 2$. The Mallows' distance $d_M(x_{1},x_{2})$ between $x_1$ and $x_2$ is defined by 
    \begin{equation}
    \label{Dist1:Wasser_q}
        d_{M}(x_{1},x_{2})= \left(\int_{0}^{1}{\left(F_1^{-1}(t) - F_2^{-1}(t)\right)^2}dt\right)^{1/2}.
    \end{equation}
\end{definition}

Note that \eqref{Dist1:Wasser_q} is a distance between the quantile functions of the microdata within the macrodata associated to $x_1$ and $x_2$. In this sense, two symbolic intervals are at Mallows' distance $0$ if and only if they have the same quantile function. In fact, using the model described in Definition~\ref{Def:Model_Micro_Macro}, we can show that the quantile function of the microdata is a transformation of the quantile function of the latent random variable defined in $[-1,1]$. Firstly, we prove the following auxiliary lemma.

\begin{lemma}
\label{Prp:quantileF}
Let $W$ be a real-valued random variable, let $Z=a+b~W$, with $a\in\R$, $b\in \R^+_0$, and let $F_W$ and $F_W^{-1}$ be the distribution function and quantile function of $W$, respectively. Then, the quantile function of $Z$ is 
$F_Z^{-1}(t) = a + b~F_W^{-1}(t),\ t \in (0,1].$
\end{lemma}
\begin{proof}
See Appendix~\ref{AppendixA1}.
\end{proof}

Let $\Tilde{V}_i = V_{i}|\left(C_i = c_i, R_i = r_i\right) = c_i + {r_i}~\Uc_i/2$ be the random variable that describes the microdata within the macrodata of $x_i$, where $\Uc_i = U_{i}|\left(C_i = c_i, R_i = r_i\right)$, $i=1,2$. Recall that under the assumption of independence between $U_i$ and $(C_i, R_i)^T$, $\Uc_i$ and $U_i$ are identically distributed. We therefore use $U_i$ to simplify the notation. Since, for all $i$, $r_i \geq 0$, we can use the Lemma~\ref{Prp:quantileF} to show that 
\begin{equation}\label{transformquantile}
F_i^{-1}(t) = c_i + \frac{r_i}{2}~F_{\U_i}^{-1}(t),\; i=1,2,
\end{equation}
where $F_i^{-1}$ and $F_{\U_i}^{-1}$ are the quantile functions of $\Tilde{V}_i$ and $\U_i$, respectively.

Using $\eqref{transformquantile}$, we can compute the Mallows' distance between $x_1$ and $x_2$ as a function of the centres and the ranges of the two symbolic objects, and some quantities related to the latent random variables  $\U_1$ and $\U_2$.
  \begin{theorem}
    \label{Th:Dist_Wasser2_Geral}
     Let ${x}_1,\,{x}_2$ be symbolic intervals such that $x_i=\left(c_i,r_i, F_{{U}_i}\right)$, $i=1,2$, where $\U_i$ is an absolutely continuous random variable with support $[-1,1]$, assumed to have finite second moment.
     Then, the square of the Mallows' distance between $x_1$ and $x_2$ is
    \begin{align}
        d_M(x_{1},x_{2})^2\ &=\  (c_1-c_2)^2  + (c_1-c_2) \big( r_1\E{(\U_1)} - r_2\E{(\U_2)} \big) \notag\\
        &+\ \frac{ r_1^2}{ 4}\E{(\U_1^2)}\ 
        +\ \frac{ r_2^2}{ 4}\E{(\U_2^2)}
        \ -\  \frac{ r_1 r_2}{ 2}~{\mathcal E}{(\U_1 , \U_2)},\label{Eq:Wass_general_moments}
    \end{align}   
where ${\mathcal E}(\U_1,\U_2) = \int_0^1 {F_{\U_1}^{-1}(t)F_{\U_2}^{-1}(t)}\ dt$. 
\end{theorem}

\begin{proof}
See Appendix~\ref{AppendixA2}.
\end{proof}

Observe that the quantity ${\mathcal E}(\U_1,\U_2)$ is not the usual $\E{(\U_1 \U_2)}$. However, it can be seen as a cross-moment of order $2$ that depends on the quantile functions of $\U_1$ and $\U_2$, since
${\mathcal E}(\U_1,\U_2) = \int_0^1 {F_{\U_1}^{-1}(t)F_{\U_2}^{-1}(t)}\ dt = \E{\left(F_{\U_1}^{-1}(T)F_{\U_2}^{-1}(T)\right)}$
where $T$ is a real-valued random variable following a continuous uniform distribution defined in $[0,1]$.
Moreover, if $\U_1$ and $\U_2$ are identically distributed then
${\mathcal E}(\U_1,\U_2)=\E(\U_1^2)=\E(\U_2^2)$.

\begin{corollary}\label{Cor:Mallows_linked_IrpinoVerde}
    Under the conditions of Theorem~\ref{Th:Dist_Wasser2_Geral}
\begin{equation}
d_M(x_{1},x_{2})^2\ =\ (\mu_1-\mu_2)^2 + (\sigma_1-\sigma_2)^2 + 2\sigma_1\sigma_2(1-\rho_{12}),\label{Eq:Wass_mu_sigmas}
   \end{equation} 
where $\mu_i = \E{(\tilde{V}_i)} = c_i +  r_i\,\E{({U}_i)}/2$,  \hspace{2pt} $\sigma_i^2 =  r_i^2\,\Var{(\U_i)}/4$,
and 
\begin{equation*}
\rho_{12} = \frac{\int_0^1{F_1^{-1}(t)F_2^{-1}(t)}dt - \mu_1\mu_2}{\sigma_1\sigma_2} = \frac{{\mathcal E}{(\U_1 , \U_2)}-\E{(\U_1)}\E{(\U_2)}}{\sqrt{\Var{(\U_1)}\Var{(\U_2)}}}.
\end{equation*}
\end{corollary}
\begin{proof}
    See Appendix~\ref{AppendixA2.B}. 
\end{proof}

In \cite[Prop. 2]{Irpino.verde:2015} the authors deduced the formulation \eqref{Eq:Wass_mu_sigmas} for the Mallows' distance and identified $\rho_{12}$ as the correlation coefficient between two quantile functions. This can be confirmed by considering 
$\E{(\U_i)}=\E{\left(F_{\U_i}^{-1}(T)\right)}$, $\E{(\U_i^2)}=\E{\left(\left(F_{\U_i}^{-1}(T)\right)^2\right)}$ and ${\mathcal E}(\U_1 , \U_2)=\E{\left( {F_{\U_1}^{-1}(T)F_{\U_2}^{-1}(T)}\right)}$, $T \sim \rm{Unif}(0,1)$.

\begin{remark}
The Mallows' distance between a symbolic interval $x_1$ and a point $x_2$ is given by $d_M(x_{1},x_{2})^2\ =\  (c_1-c_2)^2\  +\ (c_1-c_2)\  r_1 \E{(\U_1)}\ +\  r_1^2\E{(\U_1^2)}/4$. This follows from \eqref{L2quantile} and the fact that $r_2 = 0$ and $F_{\U_2}^{-1}(t) = 0,\ t \in (0,1]$, according to Definition~\ref{Def:Model_Micro_Macro}. This expression can also be obtained from \eqref{Eq:Wass_general_moments} by setting $r_2=0$. 
\end{remark}

The main challenge in Theorem~\ref{Th:Dist_Wasser2_Geral} is to find the quantity ${\mathcal E}{(\U_1 , \U_2)}$ by calculating the integral of the product of the quantile functions of $\U_1$ and $\U_2$. In the next example, we compute the square of the Mallows' distance between two symbolic intervals, where the distribution of the latent random variables $\U_1$ and $\U_2$ is known.

\begin{example}
     Let ${x}_1$ and ${x}_2$ be symbolic intervals such that $x_i=\left(c_i,r_i, F_{\U_i}\right)$, $i=1,2$, where $\U_1$ follows a continuous uniform distribution $\textnormal{Unif}(-1,1)$ and $\U_2$ follows a symmetric triangular distribution $\textnormal{Triang}(-1,1,0)$ (see  Figures \hyperref[Fig:CUnif]{2b} and \hyperref[Fig:InvT]{2c}).
     Since both distributions are symmetric, $\E(\U_1) = \E(\U_2) = 0$. Furthermore, it can be easily shown that $\E(\U_1^2) = 1/3$ and $\E(\U_2^2) = 1/6$. It only remains to compute the quantity ${\mathcal E}{(\U_1 , \U_2)}$. Noting that the quantile function of $\U_1$ is $F_{\U_1}^{-1}(t) = 2t-1$, $t \in (0,1]$, the quantile function of $\U_2$ is $F_{\U_2}^{-1}(t) = -1 + \sqrt{2t}$, if $t \in (0,1/2]$, and $F_{\U_2}^{-1}(t) = 1 - \sqrt{2(1-t)}$, if $t \in (1/2,1]$, we have
\begin{align*}
    {\mathcal E}{(\U_1 , \U_2)}\ &=\ \int_0^1 {F_{\U_1}^{-1}(t)F_{\U_2}^{-1}(t)}\ dt\\ 
    &=\ \int_0^{\frac{1}{2}} {\left(2t-1\right)\left(-1+\sqrt{2t}\right)}\ dt\ +\ \int_{\frac{1}{2}}^1 {\left(2t-1\right)\left(1-\sqrt{2(1-t)}\right)}\ dt\ =\ \frac{7}{30}.
\end{align*}
Hence, the square of the Mallows' distance between $x_1$ and $x_2$ is
\begin{equation*}
    d_M(x_1,x_2)^2\ =\ (c_1-c_2)^2\ 
        +\ \frac{1}{12}r_1^2\ +\ \frac{1}{24}r_2^2\ -\ \frac{7}{60}r_1r_2.
\end{equation*}

\end{example}

When $\U_1$ and $\U_2$ are identically distributed, ${\mathcal E}{(\U_1 , \U_2)}=\E(\U_1^2)=\E(\U_2^2)$. 
This leads to interesting simplifications of the results stated in Theorem~\ref{Th:Dist_Wasser2_Geral}, as shown in the next corollary.
 \begin{corollary}
    \label{Crll:Dist_Wasser2}
     Under the conditions of Theorem~\ref{Th:Dist_Wasser2_Geral}, when $\U_1$ and $\U_2$ are identically distributed, 
     the square of the Mallows' distance between $x_1$ and $x_2$ is
    \begin{align}
    d_M(x_1,x_2)^2\
    =\ (c_1-c_2)^2 + \frac{\E{(\U_1^2)}}{4}~(r_1-r_2)^2 + \E{(\U_1)}~(c_1-c_2)\left(  r_1 - r_2 \right).\label{Eq:D2WasserE(U)qq}
    \end{align} 
 Additionally, if $\U_1$ and $\U_2$ are also symmetric random variables, then $\E{(\U_1)} = \E{(\U_2)} = 0$ and 
     \begin{equation}
    d_M(x_1,x_2)^2\ =\ (c_1-c_2)^2+ \delta (r_1-r_2)^2,  \label{Eq:DW_SymIID}   
    \end{equation} 
 where $\delta =  \Var{(\U_1)}/4$.
  
\end{corollary}

\begin{remark}
    It is possible to establish that the parameter $\delta = \Var{(\U_1)}/4$ (and also $\E{(\U_1^2)}/4$) takes values in the interval $[0,1/4]$. Trivially, the variance is a non-negative number. To show the upper bound, let $f_{\U_1}$ be the probability density function of $\U_1$, then %
    \begin{equation*}
    \Var{(\U_1)}\ \leq \ \E(\U_1^2)\ = \int_{-1}^{1}\!x^2f_{\U_1}(x)\ dx\ \leq \int_{-1}^{1}f_{\U_1}(x)\ dx  \ = 1.
\end{equation*}

It is also worth noting that in the formulation \eqref{Eq:DW_SymIID} (and \eqref{Eq:D2WasserE(U)qq}) of the Mallows' distance, the weight associated with the squared distance between the ranges ($\delta=\Var{(\U_1)}/4$ in \eqref{Eq:DW_SymIID} and $\E{(\U_1^2)}/4$ in \eqref{Eq:D2WasserE(U)qq}) is always smaller than or equal to 1/4. 
This emphasises the lesser role of the ranges when compared to the centres. 
If the variance of $\U_1$ goes to zero, then the microdata is more and more concentrated around the centre of the interval $c_1$. 
\end{remark}

To have a geometric interpretation of the Mallows' distance \eqref{Eq:DW_SymIID} (symmetric and identically distributed latent random variables), we consider the sets ${\mathcal A}_\delta$ of symbolic intervals whose distance to the symbolic interval $x_0=([-3,5],F_{\U_1})$ is one unit, indexed by the parameter $\delta \in [0,1/4]$ used in the distance. Since the centre and range of $x_0$ are, respectively, $c_0=1$ and $r_0=8$, we have, for $\delta=\Var{(\U_1)}/4$,
\begin{equation*}
{\mathcal A}_\delta=\left\{x=(c,r,F_{\U_1}) :\, d_M(x, x_0) = 1\right\}=\left\{(c,r,F_{\U_1}) :\, (c-1)^2 + \delta (r-8)^2 = 1\right\}.
\end{equation*}

\begin{figure}[h]
\centering
\includegraphics[alt = {Ellipses representing the set of points at Mallows' distance one. For symmetric distributions with smaller variance the ellipses become more elongated.}, width=0.65\textwidth]{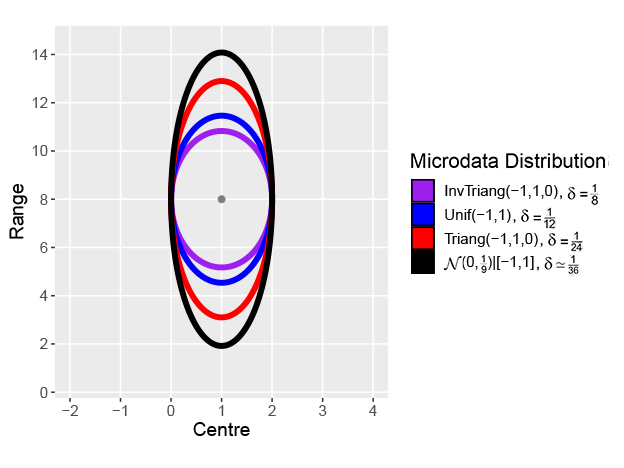}
\caption{Set of intervals whose Mallows' distance to $x_0=(1,8,F_{\U_1})$ is one unit, according to the assumed symmetric distribution for $\U_1$}
  \label{fig:Mallows_Ellipses}
\end{figure}

Note that each of these sets describes an ellipse whose elongation is controlled by the parameter $\delta$. Figure~\ref{fig:Mallows_Ellipses} illustrates some of these ellipses for the different distributions presented in Figure~\ref{fig:distributions} and their respective value of $\delta$. The lower the variance of $\U_1$, the greater the concentration of microdata around the centres of the intervals, and the larger the area of the region whose symbolic intervals $x$ verify $d_M(x, x_0) \leq 1$.

In \cite{IrpinoRomano2007,Irpino.verde:2015} the authors derived \eqref{Eq:Wass_mu_sigmas}  based on the quantile functions related to the intervals $x_j,\ {j=1,2}$, and considered the special case when the microdata within $x_j$ followed a continuous uniform distribution,  resulting in $\delta=1/12$. In \cite{Dias.Brito:2017} a similar result was obtained for a symmetric triangular distribution leading to $\delta=1/24$ in \eqref{Eq:DW_SymIID}.

\subsection{Multivariate Mallows' distance}\label{Sec:WasserstIRp}

In SDA, the generalisation of the Mallows' distance to higher dimensions is defined using the Mallows' distance between each component of the two vectors.

\begin{definition}\label{Def:Dist_p_Mink}
    Let $\bm{x}_1=(x_{11},\ldots,x_{1p})$ and $\bm{x}_2=(x_{21},\ldots,{x}_{2p})$ be symbolic hyperrectangles.
    In addition, let $d_M(x_{1i},x_{2i})$, ${i=1,\ldots,p}$, be the univariate Mallows' distance between the symbolic intervals $x_{1i}$  and $x_{2i}$. The Mallows' distance between $\bm{x}_1$ and $\bm{x}_2$ is
    \begin{equation}
    \label{Dist:Minkowski}
        d_M(\bm{x}_1,\bm{x}_2)=\left(\sum_{i=1}^p d_M(x_{1i},x_{2i})^2\right)^{1/2}.
    \end{equation}
\end{definition}

Since the Mallows' distance between symbolic hyperrectangles is defined using the distance between each component, we are interested in comparing the componentwise quantile functions of the microdata. As such, a symbolic hyperrectangle is uniquely identified by the vector of the centres, the vector of the ranges, and the collection of quantile functions of the microdata in each dimension. We use the notation ${\bm{x}_j = \left(\bm{c}_j, \bm{r}_j, F_{\U_{j1}},\ldots,F_{\U_{jp}}\right)}$, where each component is of the form $x_{ji} = \left(c_{ji}, r_{ji}, F_{\U_{ji}}\right),\  j=1,2$,\ $i = 1,\ldots,p$.

The following results require that, within the same dimension, the latent random variables of the intervals are identically distributed. In this scenario, we introduce a latent random variable, $\U_i,\ i=1,\ldots,p,$ whose distribution is the same as $\U_{ji}$, $j=1,2$, such that $\E(\U_{1i}) = \E(\U_{2i}) = \E(\U_i)$ and $\E(\U_{1i}^2) = \E(\U_{2i}^2) = \E(\U_i^2)$. We write $\bm{x}_j = \left( \bm{c}_j, \bm{r}_j, F_{\U_{1}},\ldots,F_{\U_{p}} \right)$, or simply $\bm{x}_j = (\bm{c}_j, \bm{r}_j, F_{\bm{\U}})$, $j=1,2$, if no confusion arises. Given a vector $\bm{v}$, we define $\diag(\bm{v})$ to be the diagonal matrix whose main diagonal is $\bm{v}$. In addition, given a matrix $\bm{A}$, we define $ \diag(\bm{A})$ as the diagonal matrix with the same main diagonal as $\bm{A}$.
\begin{theorem}
\label{Th:Dist_Wasser2_id}
For $j=1,2$ let $\bm{x}_j = \left(\bm{c}_j, \bm{r}_j, F_{\bm{\U}}\right)$, with $\bm{c}_j = (c_{j1}, \ldots, c_{jp})^T \in \R^p$, 
$\bm{r}_j = (r_{j1}, \ldots, r_{jp})^T  \in (\R^+_0)^p$, and $\bm{\U}=(\U_1,\ldots,\U_p)^T$ a random vector with support $[-1,1]^p$. 
Then, the square of the Mallows' distance between $\bm{x}_1$ and $\bm{x}_2$ is
\begin{align}
      d_M(\bm{x}_{1},\bm{x}_{2})^2
        \ &=\ \sum_{i=1}^p \left((c_{1i}-c_{2i})^2          
        + \frac{\E{(\U_i^2)}}{4}  (r_{1i}-r_{2i})^2
        + \E{(\U_i)} (c_{1i}-c_{2i})(r_{1i} - r_{2i})\right) \notag\\  
        &=\ (\bm{c}_1-\bm{c}_2)^T(\bm{c}_1-\bm{c}_2)
        + \left( \bm{r}_{1}-\bm{r}_{2}\right)^T \bm{\Delta} \left( \bm{r}_{1}-\bm{r}_{2}\right)
        +(\bm{c}_1-\bm{c}_2)^T \bm{\Psi} (\bm{r}_1-\bm{r}_2), \label{Eq:MD2WasserE(U)qq}
\end{align}
where $\bm{\Delta}=\diag(\delta_1,\ldots,\delta_p)$ with $\delta_i=\E(\U_i^2)/4$, and
$\bm{\Psi}=\diag\big(\E(\U_1),\ldots,\E(\U_p)\big)$.
\end{theorem}
\begin{proof}
    The result follows from replacing the general expression of Corollary~\ref{Crll:Dist_Wasser2} in each component of the Mallows' distance \eqref{Dist:Minkowski} and using matrix notation.
\end{proof}

\begin{remark}
    Observe that, when the random variables $\U_i$ are symmetric, $\bm{\Psi}=\bm{0}$ and \eqref{Eq:MD2WasserE(U)qq} becomes
\begin{equation*}
    d_M(\bm{x}_{1},\bm{x}_{2})^2\ =\ (\bm{c}_{1}-\bm{c}_{2})^T(\bm{c}_{1}-\bm{c}_{2}) + \left( \bm{r}_{1}-\bm{r}_{2}\right)^T \bm{\Delta} \left( \bm{r}_{1}-\bm{r}_{2}\right),
\end{equation*}
where $\delta_i=\Var(\U_i)/4$. In particular, if $\delta_i=\delta$, $i=1,\ldots,p$, we obtain $\bm{\Delta}=\delta \bm{I}_p$, where $\bm{I}_p$ represents the  identity matrix of dimension $p$, and

\begin{equation*}
    d_M(\bm{x}_{1},\bm{x}_{2})^2\ =\ (\bm{c}_{1}-\bm{c}_{2})^T(\bm{c}_{1}-\bm{c}_{2}) + \delta \left( \bm{r}_{1}-\bm{r}_{2}\right)^T \left( \bm{r}_{1}-\bm{r}_{2}\right).\label{Eq:L2Wasser_Multi_oneDelta}
\end{equation*}

\end{remark}

As an alternative to \eqref{Eq:MD2WasserE(U)qq}, Corollary~\ref{corollary:mahalob} proposes the identification of the Mallows' distance $d_M$ with a weighted Euclidean distance $d_{H}$ in the space of the joint vector of the centres and ranges. This follows from the fact that within each dimension the latent random variables have the same distribution function, yielding that the only difference between the hyperrectangles is in the centres and the ranges.

\begin{corollary}\label{corollary:mahalob}
    Under the conditions of Theorem~\ref{Th:Dist_Wasser2_id}, the square of the Mallows' distance \eqref{Eq:MD2WasserE(U)qq} between $\bm{x}_1 = \left(\bm{c}_1, \bm{r}_1, F_{\bm{\U}}\right)$ and $\bm{x}_2 = \left(\bm{c}_2, \bm{r}_2, F_{\bm{\U}}\right)$, can be expressed as 
\begin{equation} \label{Eq:Mallows_Mahalob}
      d_M(\bm{x}_{1},\bm{x}_{2})^2
        = d_{ H}(\bm{y}_{1},\bm{y}_{2})^2 = \left(\bm{y}_1-\bm{y}_2\right)^T \bm{H} (\bm{y}_1-\bm{y}_2),\ \text{where} \ \bm{H}=\left(
\begin{array}{cc}
\bm{I}_p & \frac{1}{2}\bm{\Psi}\\
\frac{1}{2}\bm{\Psi} & \bm{\Delta}\\
\end{array}
\right)\ ,
\end{equation}

\noindent and $\bm{y}_i = (\bm{c}_i^T, \bm{r}_i^T)^T$ represents the macrodata of $\bm{x}_i$.

If $\bm{H}$ is positive definite (i.e., the random variables ${\U}_i$, $i=1,\ldots,p$, are non-degenerate), 
then $d_{H}(\bm{y}_{1},\bm{y}_{2})$ is a weighted Euclidean distance in $\R^{p}\times (\R_0^+)^p$, and $\bm{H}^{-1}$ can be interpreted as a covariance matrix where 
 \begin{equation*}
    \bm{H}^{-1}=\left(
\begin{array}{cc}
\bm{I}_p + \frac{1}{4}\bm{\Psi}^2\bm{Q} & -\frac{1}{2}\bm{\Psi}\bm{Q}\\
-\frac{1}{2}\bm{\Psi}\bm{Q} & \bm{Q}\\
\end{array}
\right),  
 \end{equation*}
and $\bm{Q}=(\bm{\Delta}-\frac{1}{4} \bm{\Psi}^2)^{-1} = 4\diag(\Var(\U_1)^{-1},\ldots,\Var(\U_p)^{-1})$, $\Var(\U_i)>0$, $i = 1,\ldots,p$.
\end{corollary}
    
\begin{proof}
    See Appendix~\ref{AppendixA3}.
\end{proof}

\begin{remark}
    If any ${\U}_i$ is a degenerate random variable —- without loss of generality, say ${\U}_1$ for simplicity, then $X_1$ is a conventional variable. In that case, we define a weighted Euclidean distance in a space of reduced dimension, $\R^{p}\times (\R_0^+)^{p-1}$, by removing the $(p+1)$-th row and column of $\bm{H}$.
\end{remark}

\section{Location, scale, and association between interval-valued variables }\label{Sec4:moments}

\cite{Irpino.verde:2015} proposed an approach to derive interval and histogram-valued descriptive measures for location, scale, and association between two symbolic characteristics measured on the same set of objects.  The authors relied on the (sample) Fréchet mean, also known as the (sample) barycentre. Given a set of points in a metric space, the Fréchet mean is the point that minimises the weighted sum of the squares of the distance to all given points \citep{2011:Agueh.Carlier}. The minimum of this sum is called the Fréchet variance. The definition of Fréchet mean can be extended to the population if, instead of the weighted sum, we consider the expected value. In the case of an interval-valued random vector, $\bm{X}$, we can define the population Fréchet mean or population barycentre, using the Mallows' distance, $d_M$. Given a matrix $\bm A$, we define $\tr(\bm{A})$ as its trace, i.e., the sum of the elements on the main diagonal.

\begin{definition}
\label{Eq:PopBaryc0}
    Let $\bm{X}$ be an interval-valued random vector characterised by  $(\bm{C},\bm{R},F_{\bm{U}})$, where $\bm{C}$ and $\bm{R}$ are assumed to have finite expected values, $\bm{\mu}_C$ and $\bm{\mu}_R$,  covariance matrices, $\bm{\Sigma}_{CC}$ and $\bm{\Sigma}_{RR}$, respectively, and  $\bm{\Sigma}_{CR}$ be the covariance matrix between $\bm{C}$ and $\bm{R}$. Let ${\bm{U}=(U_1,\ldots,U_p)^T}$ be a real-valued random vector of independent random variables with distribution function $F_{\bm{U}}$, where $\bm{U}$ is assumed to be independent from $(\bm{C}, \bm{R})$, let $\bm{\Delta}=\diag(\delta_1,\ldots,\delta_p)$, where ${\delta_i=\E{(U_i^2)}/4}$, $i=1,\ldots,p$, and let 
    $\bm{\Psi}=\diag\left(\E(U_1),\ldots,\E(U_p)\right)$. 
    The population barycentre of $\bm{X}$, denoted by $\bm{\mu}_B$ and characterized by $(\bm{c},\bm{r},F_{\bm{U}})$ is the interval object that globally minimises the function 
    \begin{eqnarray}
        g(\bm{x}) &=& g(\bm{c},\bm{r};F_{\bm{U}}) = \E{\left(d_M(\bm{X},\bm{x})^2\right)}\ \nonumber\\
        &=& \E\left( (\bm{C}-\bm{c})^T(\bm{C}-\bm{c}) 
        \ +\  (\bm{R}- \bm{r})^T\bm{\Delta}(\bm{R}-\bm{r}) 
        \ +\ (\bm{C}-\bm{c})^T \bm{\Psi} (\bm{R}-\bm{r})\right). \label{Eq:Barycenter.Pop}
    \end{eqnarray}
\end{definition}

The solution of the minimisation of \eqref{Eq:Barycenter.Pop} is given in the following theorem.

\begin{theorem}
\label{Th:Barycenter.Pop}
Under the conditions of Definition~\ref{Eq:PopBaryc0}, the population barycentre of $\bm{X}$ is 
\begin{equation}
\label{Eq:PopBaryc}
    \bm{\mu}_B\ =\ (\bm{\mu}_C,\bm{\mu}_R, F_{\bm{U}})
\end{equation}  
and the corresponding Fréchet variance is 
\begin{equation}\label{Eq:FrechetVar_Multivar}
   V_F(\bm{\mu}_B)\ =\ \E\left(d_M(\bm{X},\bm{\mu}_B)^2\right) 
 \ =\ \tr(\bm{\Sigma}_{CC}
 +\bm{\Delta} \bm{\Sigma}_{RR}
 +\bm{\Sigma}_{CR}\bm{\Psi}).
  \end{equation}
\end{theorem}
\begin{proof}
    See \ref{Prf:Barycenter.Pop}.
\end{proof}

As a particular case, notice that, if we observe a sample of $\bm{X}$, say $\bm{x}_1,\ldots,\bm{x}_n$, 
Theorem~\ref{Th:Barycenter.Pop} guarantees that the (sample) barycentre $\bar{\bm{x}}_B =
    (\bm{c}_B, \bm{r}_B,F_{\bm{U}})$ verifies: 
\begin{align}\label{Eq:SampleBaryOF}
    \argmin_{(\bm{c},\bm{r})}\, \frac{1}{n}\sum_{j=1}^n \left[(\bm{c}_j-\bm{c})^T(\bm{c}_j-\bm{c}) 
     \ +\ (\bm{r}_j - \bm{r})\bm{\Delta}(\bm{r}_j-\bm{r})
     \ +\ (\bm{c}_j-\bm{c}) \bm{\Psi} (\bm{r}_j-\bm{r}) \right],
\end{align}
     leading to $\bm{c}_B=
     \overline{\bm{c}}_n = \sum_{j=1}^{n}\bm{c}_{\,j}/n$ and 
     $\bm{r}_B= \overline{\bm{r}}_n = \sum_{j=1}^{n}\bm{r}_j/n$. In short, $\bar{\bm{x}}_B\ = (\overline{\bm{c}}_n, \overline{\bm{r}}_n,F_{\bm{U}})$.

\cite{Irpino.verde:2015} propose the definition of sample variance and sample covariance between two interval-valued variables based on the Mallows' distance that follows, without any link to the model presented in Definition~\ref{Def:Model_Micro_Macro}.
\begin{definition}
    \label{Def:SampleCov_B}
    Let $({x}_{11},x_{12})^T, \ldots,({x}_{n1},x_{n2})^T$  be a sample of size $n$ from the bivariate interval-valued random vector $(X_1,X_2)^T$,  with microdata distribution $F_{ji}$  having  support on the macrodata of ${x}_{ji}$, for $j=1,\ldots,n$ and $i=1,2$.
    The associated sample barycentre is $(\overline{\bm{c}}_{1n},\overline{\bm{c}}_{2n}, \overline{\bm{r}}_{1n},\overline{\bm{r}}_{2n} ,  F_{B_1},F_{B_2})$, where $F_{B_i}$ is the distribution function of the microdata  with support $[\overline{\bm{c}}_{in}-\overline{\bm{r}}_{in}/2, \overline{\bm{c}}_{in}+\overline{\bm{r}}_{in}/2 ]$, for $i=1,2$.
    The sample covariance between $X_1$ and $X_2$ is 
\begin{equation}\label{Eq:SampleCov_B}
    s_{12,B}\ =\ \widehat{{\rm Cov}}_B{(X_1,X_2)}\ =\ \frac{1}{n} \sum_{j=1}^{n}\, \int_0^1 \left( F_{j1}^{-1}(t)-F_{B_1}^{-1}(t)\right) \left(F_{j2}^{-1}(t)-F_{B_2}^{-1}(t) \right) dt,
\end{equation}
Additionally, the sample variance of $X_i$ is $s_{ii,B}=\widehat{{\rm Cov}}_B{(X_i,X_i)}$. 
\end{definition}

In the particular case $p=1$,  $s_{11,B}$ is the minimum value of the objective function in \eqref{Eq:SampleBaryOF} and corresponds to the univariate sample Fréchet variance.

The population counterpart to the symbolic sample covariance can be defined by a natural adaptation of \eqref{Eq:SampleCov_B}.

\begin{definition}
    \label{Def:Cov_B}
    Let $X_1 = (C_1,R_1,F_{U_1})$ and $X_2 = (C_2,R_2,F_{U_2})$ be two interval-valued random variables with baricentres $\mu_{B_1}$ and $\mu_{B_2}$, respectively. Let $G_{i}^{-1}(t) = C_i + R_iF_{U_i}^{-1}(t)/2$, $t\in[0,1]$, be the random variable whose realisations on specific intervals are the quantile functions of the microdata within, where $F_{U_i}^{-1}$ is the quantile function of the latent random variable $U_i$. Let $F_{B_i}^{-1}(t)$, $t\in[0,1]$, be the quantile function of the microdata in the real interval characterizing $\mu_{B_i}$, $i=1,2$. The covariance ${\rm Cov}_B{(X_1,X_2)}$ between $X_1$ and $X_2$ is defined by
    \begin{equation}
        {\rm Cov}_B{(X_1,X_2)}\ = \ \E\left(  \int_0^1 \left( G_{1}^{-1}(t)-F_{B_1}^{-1}(t)\right) \left(G_{2}^{-1}(t)-F_{B_2}^{-1}(t) \right) dt \right).
    \end{equation}
Furthermore, the variance of $X_i$ is ${\rm Cov}_B{(X_i,X_i)}$, $i=1,2$. 
\end{definition}

Considering that $\Tilde{V}_{B_i}=\mu_{C_i}+\U_{B_i}\hspace{2pt}\mu_{R_i}/{2}$ (see Definition~\ref{Def:Model_Micro_Macro}) and noting that $U_i$ and $\U_{B_i}$ are identically distributed, we can simplify the bivariate covariance matrix.

\begin{corollary}\label{Crll:PopCov}
    Under the conditions of Definition~\ref{Eq:PopBaryc0}, the covariance between two interval-valued random variables $X_1 = (C_1,R_1,F_{U_1})$ and $X_2 = (C_2,R_2,F_{U_2})$ with barycentres $\mu_{B_1}$ and $\mu_{B_2}$, respectively, is 
\begin{align*}   
    {\rm Cov}_B{(X_1,X_2)}\ &= \ \E\left(  \int_0^1 \left( G_{1}^{-1}(t)-F_{B_1}^{-1}(t)\right) \left(G_2^{-1}(t)-F_{B_2}^{-1}(t) \right) dt \right)\\
    &=\ \Cov(C_1,C_2)
    \ +\ \frac{{\mathcal E}(U_1 , U_2)}{4}\Cov(R_1,R_2)
    \\ &+\ \frac{\E(U_2)}{2}\Cov(C_1,R_2)\ +\ \frac{\E(U_1)}{2}\Cov(C_2,R_1),
\end{align*}
where ${\mathcal E}(U_1 ,U_2)=\int_0^1 F_{U_1}^{-1}(t) F_{U_2}^{-1}(t)\ dt$.
Moreover, 
\begin{equation} 
\label{Eq:PopVar}
    {\rm Var}_B{(X_1)}\ =\ {\rm Cov}_B{(X_1,X_1)}
    \ =\ \Var(C_1)
    \ +\ \frac{\E(U_1^2)}{4}\Var(R_1)\ +\ \E(U_1)\hspace{1pt}\Cov(C_1,R_1).
\end{equation}
    
\end{corollary}
\begin{proof}
    See \ref{Prf:PopCov}.
\end{proof}

Applying this result to a $p$-dimensional random vector results in a covariance matrix, ${{\rm Var}_B{(\bm{X})}=\bm{\Sigma}_B}$, as stated in the following corollary. Here, we introduce the notation $[\bm{A}]_{ij}$ to represent the entry $(i,j)$ of matrix $\bm{A}$.

\begin{corollary}\label{Crll:Bary_CovMatrix}
    Under the conditions of Definition~\ref{Eq:PopBaryc0}, the covariance matrix of an $p$-dimensional interval-valued random vector $\bm{X}$ is
    \begin{equation}\label{Eq:L2_cov_matrixSchur}
     {\rm Var}_B{(\bm{X})}\ =\ \bm{\Sigma}_B\ =\ \bm{\Sigma}_{CC}
     \ +\  \frac{1}{4}\bm{\mathfrak E}_{UU}\bullet \bm{\Sigma}_{RR}
     \ +\ \frac{1}{2}\bm{\Sigma}_{CR} \bm{\Psi}
     \ +\ \frac{1}{2}\bm{\Psi}\bm{\Sigma}_{RC},
    \end{equation}  
   where $\bm{\Sigma}_{CC}$ and $\bm{\Sigma}_{RR}$  are the respective covariance matrices of $\bm{C}$ and $\bm{R}$, $\bm{\Sigma}_{CR}=\bm{\Sigma}_{RC}^T$ is the covariance matrix between $\bm{C}$ and $\bm{R}$, $\bm{\Psi}=\diag\left(\E(U_1),\ldots,\E(U_p)\right)$, ${\mathcal E}(U_i ,U_j)=\int_0^1 F_{U_i}^{-1}(t) F_{U_j}^{-1}(t)\ dt$, $[\bm{\mathfrak E}_{UU}]_{ij}={\mathcal E}(U_i , U_j)$, $i\neq j$, $[\bm{\mathfrak E}_{UU}]_{ii}={\E}(U_i^2)$, $i,j = 1,\ldots,p$,  and  
   \mbox{``\hspace{3pt}$\bullet$''} denotes the Schur (or entrywise) product of matrices. 
The corresponding correlation matrix is ${\rm Cor}_B{(\bm{X})} = \bm{D}^{-1/2}\bm{\Sigma}_B\bm{D}^{-1/2}$, 
 where $\bm{D}=\diag \left( [\bm{\Sigma}_B]_{11},\ldots,[\bm{\Sigma}_B]_{pp} \right) $.
\end{corollary}

As before, specific assumptions on the random variables $U_i$, $i = 1,\ldots,p$, lead to simpler covariance matrices $\bm{\Sigma}_B$. For example, if all $U_i$ are identically distributed to a random variable $U$, then ${{\mathcal E}(U_i , U_j)=\E(U^2)=4\delta}$, and 
\begin{equation}\label{Eq:L2_cov_matrixD}
     {\rm Var}_B{(\bm{X})}\ =\ \bm{\Sigma}_B\ =\ \bm{\Sigma}_{CC}\ +\ \delta\bm{\Sigma}_{RR}\ +\ \frac{\E(U)}{2}\left(\bm{\Sigma}_{CR}\ +\ \bm{\Sigma}_{RC}\right).
    \end{equation}  
Moreover, if $\E(U)=0$, then 
\begin{equation}\label{Eq:L2_cov_matrixD0}
     {\rm Var}_B{(\bm{X})}\ =\ \bm{\Sigma}_B\ =\ \bm{\Sigma}_{CC}\ +\  \delta\bm{\Sigma}_{RR}.
    \end{equation} 

Having in mind \eqref{Eq:FrechetVar_Multivar}, the Fréchet variance is a non-negative scaler interpreted as the total variance of the covariance matrix $\bm{\Omega}=\bm{\Sigma}_{CC} + \bm{\Delta} \bm{\Sigma}_{RR}+\bm{\Sigma}_{CR}\bm{\Psi}$. 
However, $\bm{\Omega}$ differs from the barycentre-based covariance matrix $\bm{\Sigma}_B$. Despite this distinction, they share the same trace (that is, the same total variance) because $\tr(\bm{\Sigma}_{CR} \bm{\Psi})=\tr(\bm{\Psi}\bm{\Sigma}_{RC})$ and since $[\bm{\mathfrak E}_{UU}]_{ii}={\E}(U_i^2)$, it follows that ${\tr(\bm{\mathfrak E}_{UU}\bullet \bm{\Sigma}_{RR}/4)}=\tr(\bm{\Delta} \bm{\Sigma}_{RR})$. Therefore, we conclude that $V_F(\bm{\mu}_B)=\tr(\bm{\Sigma}_B)$.

\section{Examples}
\label{Sec:Examples}
In this section, we compare several estimates of the sample mean, sample covariance, and sample correlation matrix based on three distinct datasets and consider different strategies to model the latent random variables $ U_{i}$, $i = 1,\ldots,p$.

The first example uses the \emph{credit card} dataset \citep{BillardDiday2003,BillardDiday2006,datasetCreditCard,Oliveira.et.al:2021}, where the microdata are available, and the choice of the distribution of $U_i$ was discussed in \cite{Oliveira.et.al:2021}. In this paper, the authors found evidence that these random variables followed a symmetric triangular distribution (mode zero). We revisit this problem by comparing the estimated correlation matrix of \cite{Oliveira.et.al:2021} with the one based on the barycentre approach.
The second example considers the dataset \emph{nycflights.int}, listed in \cite{datasetNewYork,DuarteS.et.al.RJournal:2021,MAINT.Data:2023}. This dataset contains information about all flights that departed from the three major New York airports in 2013. The data is aggregated by month and carrier. The microdata associated with two of the four interval variables does not suggest any obvious known family of distributions. Thus, non-parametric probability density estimators and associated quantile functions are considered.
A third dataset related to internet traffic redirection attacks \citep{Subtil.et.al:2023,Serrao.et.al:2023} is analysed under the new proposals for location and association for interval data. In this case, only microdata measures of location are available. The distribution of $ U_i$ is chosen based on the empirical knowledge of the experts and the partial information available about the microdata.

\subsection{Credit cards}
The \emph{credit card} dataset \citep{BillardDiday2003,BillardDiday2006,datasetCreditCard, Oliveira.et.al:2021} refers to five interval-valued random variables measuring the monthly expenses of three credit card users on Food ($x_1$), Social Entertainment ($x_2$), Travel ($x_3$), Gas ($x_4$), and Clothes ($x_5$), during one year, leading to a total of $n=36$ observations on $p=5$ variables.

The macrodata of the sample barycentre (see equation \eqref{Eq:PopBaryc}) is
\begin{equation*}
  \left(\,
[ 21.52,\  30.66],\;
[ 8.68 ,\  18.92],\;
[177.46 ,\ 190.47],\;
[ 20.36 ,\  29.32],\;
[ 43.37 ,\  55.26]\,\right)^T,
\end{equation*}
since \\
$\overline{\bm{c}}_n=\left(26.09,\,  13.80,\, 183.97,\,  24.84,\,  49.32\right)^T$ and 
$\overline{\bm{r}}_n=\left(9.15,  10.23,  13.01,   8.96,  11.89\right)^T$ \\
are the vectors of centres' means and ranges' means, respectively.

\begin{figure}[h]
    \centering
    \includegraphics[alt = {A pairs plot consisting of the scatter plots between the rectangles in each variable and the correlations based on two different methods. The barycentre approach and the diagonal model.}, width = 0.67\textwidth]{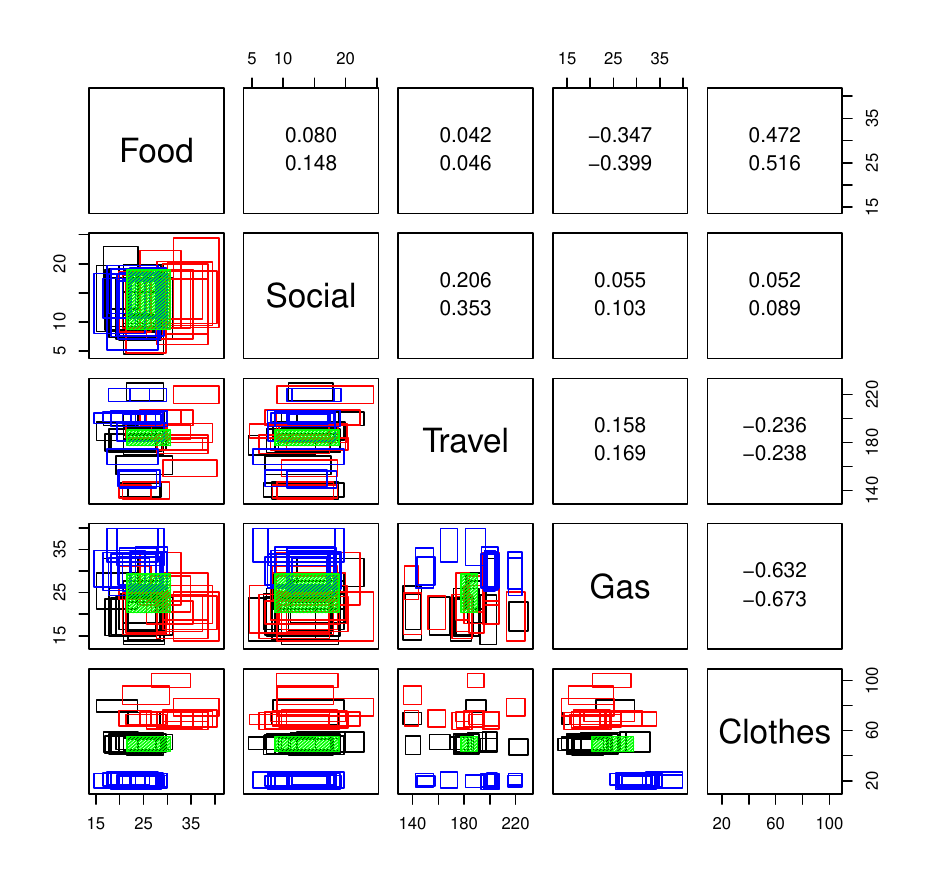}
    \caption{Symbolic bivariate scatter plots of credit card data and respective symbolic estimated correlations, based on ${\bm{\Sigma}_7=\bm{S}_{CC}+ \diag \left({\rm E}(\bm{RR}^T)\right)/24}$ (top values) and \eqref{Eq:L2_cov_matrixD0} (bottom values), assuming in both cases that $U_i \sim {\rm Triang}(-1,1,0)$. There are three subjects with monthly expenses measured over a year, coloured differently. The bivariate barycentres are in green}
    \label{fig:SymbPairsB_CCards}
\end{figure}

Figure~\ref{fig:SymbPairsB_CCards} displays a  $5\times 5$ matrix, where the entries below the main diagonal are the symbolic bivariate scatter plots of the interval-valued random variables. In green is represented the respective sample bivariate barycentres. The names of the variables appear in the main diagonal.
Figure~\ref{fig:SymbPairsB_CCards} supports the idea that the user marked in red is the one with higher expenses on Clothes ($x_5$) and Food ($x_1$). Clothes ($x_5$) is the variable that best separates the users' credit card monthly expenses. Additionally, the barycentre indicates that Travel ($x_3$) is the type of expense that users allocated the highest amount of credit card expenses, followed by Clothes ($x_5$). In opposition, Social Entertainment ($x_2$) is, on average, where the lowest amount of money is spent.
The components of the barycentre's range, $\overline{\bm{r}}_n$, are fairly similar, indicating that the inner variability among the types of expenses is also similar. The symbolic bivariate scatter plot seems to suggest a moderate positive association between Food ($x_1$) and \mbox{Clothes ($x_5$)}, a mild negative one between Gas ($x_4$) and Clothes ($x_5$), and a weaker negative association between Food ($x_1$) and Gas ($x_4$). These findings are confirmed by the associated estimated correlation values, appearing above the main diagonal of the matrix in Figure~\ref{fig:SymbPairsB_CCards}.  

The \emph{credit card} dataset was fully explored in \cite{Oliveira.et.al:2021}, where eight different symbolic estimates of covariance matrices (and respective correlation matrices)  were considered. 
Quantile-quantile plots (with 95\% pointwise envelopes) of microdata values support the assumption that the $U_i$ follow a symmetric triangular distribution, i.e., ${\rm Triang}(-1,1,m=0)$. Under the appropriate model ($k=7$ in \cite{Oliveira.et.al:2021}), it was assumed that $U_1,\ldots,U_p$ are zero mean uncorrelated random variables independent from the random vector of centres and ranges $(\bm{C}^T,\bm{R}^T)^T$. 
The estimated correlation matrix was presented in \cite[pp.~516]{Oliveira.et.al:2021} and is reproduced above the main diagonal of the matrix in Figure~\ref{fig:SymbPairsB_CCards}, in the top value of each entrance. The respective symbolic covariance matrix is then computed as $\widehat{\bm{\Sigma}}_7=\bm{S}_{CC}+ \diag \left(\widehat{{\rm E}}(\bm{RR}^T)\right)/24$, where $\widehat{{\rm E}}(\bm{RR}^T)=\bm{S}_{RR}+\overline{\bm{r}}_n\overline{\bm{r}}_n^T$, 
and $\bm{S}_{CC}$ ($\bm{S}_{RR}$) is the sample covariance matrix of the centres (ranges).

In the barycentre approach, assuming that all the $U_i$ follow a symmetric triangular distribution, we computed the estimated symbolic correlation matrix, based on $\bm{S}_B=\bm{S}_{CC}+\bm{S}_{RR}/24$ (see equation \eqref{Eq:L2_cov_matrixD0}). The values are shown above the main diagonal of the matrix of Figure~\ref{fig:SymbPairsB_CCards}, in the bottom value in each entrance. 
As expected, according to the microdata study presented in \cite{Oliveira.et.al:2021}, the two scenarios lead to a similar interpretation of the estimated correlation values.

\subsection{New York city flights interval dataset}

This example illustrates a case where the microdata are available, but the fitting of their distribution reveals an apparent major difficulty. We suggest the use of univariate kernel density estimators (KDE, see \cite{Silverman:1986}) to overcome this issue, as illustrated in this example.

\begin{figure}[!ht]
\begin{subfigure} {0.45 \textwidth} 
\includegraphics[alt = {Histogram of the latent microdata of departure delay}, width = 0.93\textwidth]{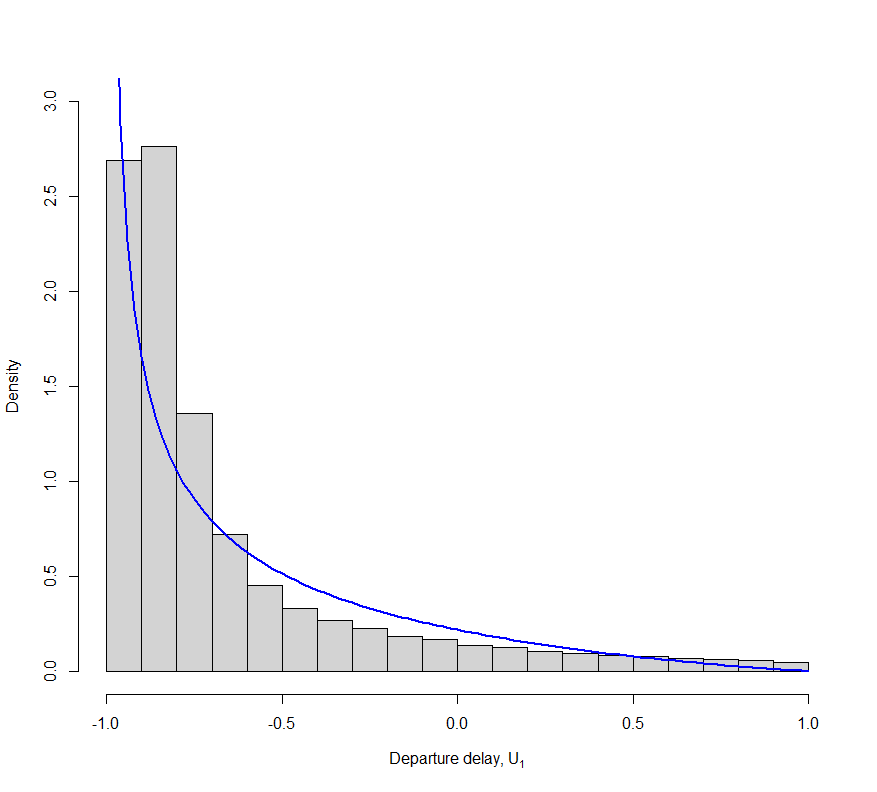} 
 \caption{Departure delay ($U_1$) and the pdf of ${(2W_1-1)}$, where $W_1 \sim {\rm Beta}(0.44, 2.15)$}
 \label{fig:NYFlights.U1}
\end{subfigure} 
\begin{subfigure}{0.45 \textwidth}
 \includegraphics[alt = {Histogram of the latent microdata of arrival delay}, width = 0.93\textwidth]{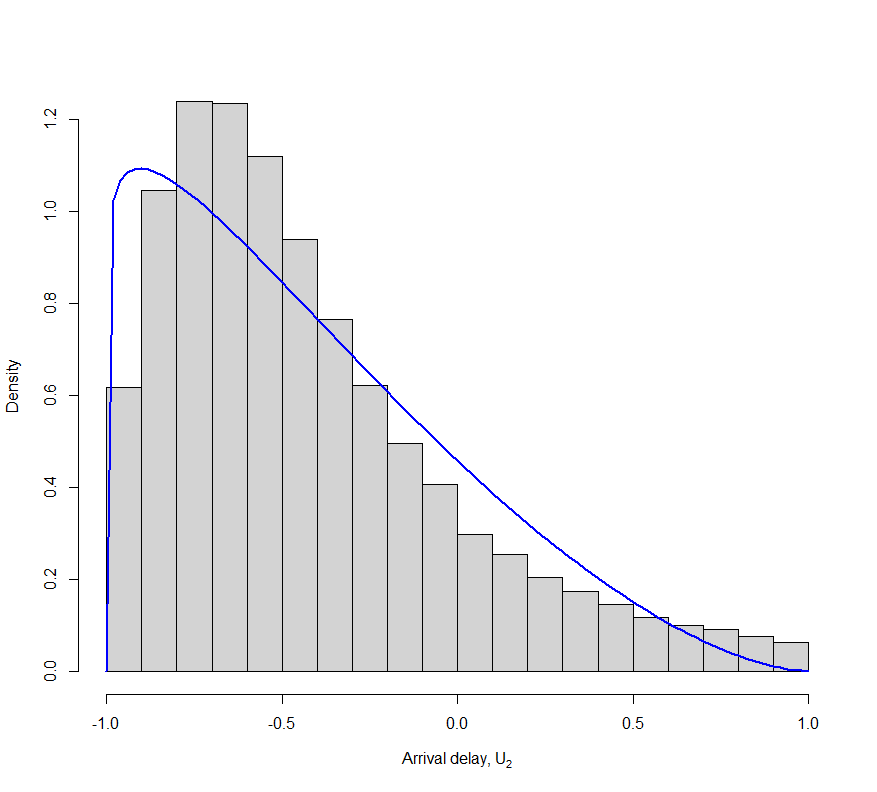}%
\caption{Arrival delay ($U_2$) and the pdf of ${(2W_2-1)}$, where $W_2 \sim {\rm Beta}(1.08, 2.65)$}
 \label{fig:NYFlights.U2}
\end{subfigure}
\hfill
\begin{subfigure} {0.45 \textwidth} 
\includegraphics[alt = {Histogram of the latent microdata of time in the air}, width = 0.93\textwidth]{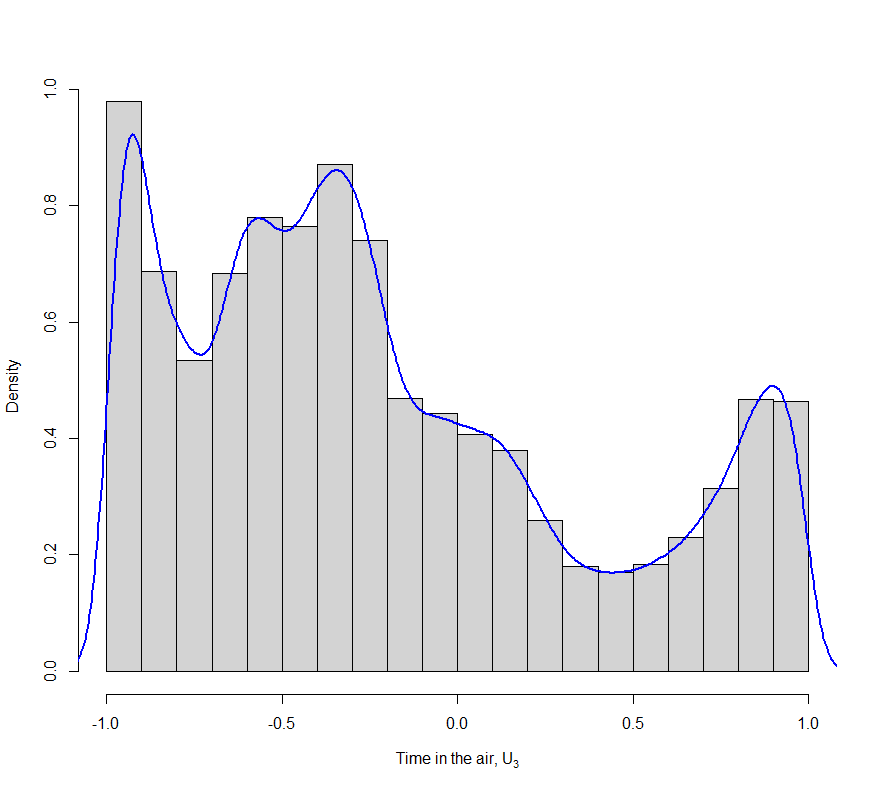} 
 \caption{Time in the air ($U_3$) and the respective KDE of $U_3$}
 \label{fig:NYFlights.U3}
\end{subfigure}
\hfill
\begin{subfigure}{0.45 \textwidth}
 \includegraphics[alt = {Histogram of the latent microdata of distance}, width = 0.93\textwidth]{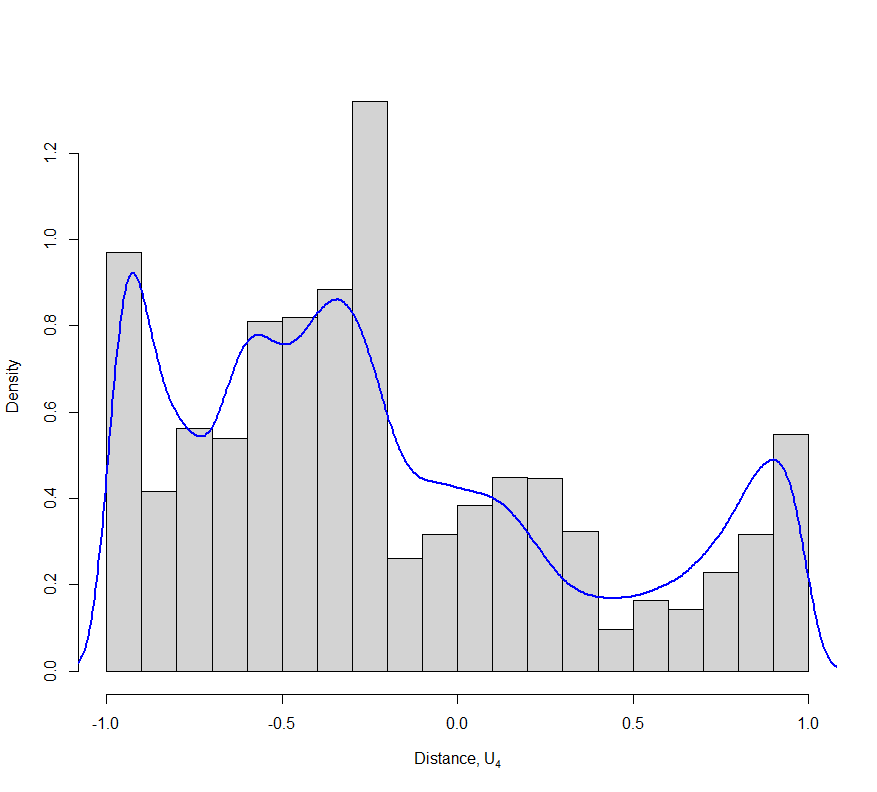}
\caption{Distance ($U_4$) and the respective KDE of $U_4$}
 \label{fig:NYFlights.U4}
\end{subfigure}
\caption{Histogram of the latent microdata related to each interval-valued variable of the New York City flights example. In the first two cases, the parameters of the Beta distributions were estimated using the method of moments. In the other two cases, kernel density estimation was used}%
 \label{fig:NYFlights_Ui}
\end{figure}

The dataset \emph{nycflights.int} \citep{datasetNewYork, DuarteS.et.al.RJournal:2021,MAINT.Data:2023} refers to all flights that departed from the three New York airports to destinations in the United States, Puerto Rico, and the American Virgin Islands, in 2013. Each flight is characterised by its departure delay ($x_1$), arrival delay ($x_2$), amount of time spent in the air ($x_3$), and distance between airports ($x_4$), for a total of $327\,345$ flights. The data were aggregated by month and carrier, leading to $n=142$ multivariate interval-valued observations and $p=4$ variables. In \cite{DuarteS.et.al.RJournal:2021,MAINT.Data:2023} authors used a robust aggregation strategy by filtering out the 5\% lowest and highest values of the microdata in each interval-valued variable. Additionally, degenerate intervals (range zero) were eliminated. The histograms of the microdata per variable are shown in Figure~\ref{fig:NYFlights_Ui}. For the first two variables, the associated latent variables are modelled as a shifted Beta distribution, i.e., $U_i=2W_i-1$, with $W_i \sim {\rm Beta}(\alpha_i, \beta_i)$, $i=1,2$. The parameters were estimated using the moment method, and the estimated probability density functions are shown in Figures Figure~\ref{fig:NYFlights.U1} and Figure~\ref{fig:NYFlights.U2}, respectively, in blue. The cases of $U_3$ and $U_4$ illustrate the difficulty in fitting the latent distribution, as shown in Figure~\ref{fig:NYFlights.U3} and Figure~\ref{fig:NYFlights.U4}. Alternatively, a univariate kernel density estimator (KDE) is used to estimate the needed quantities. The fitted KDEs were obtained using the R package \emph{kde1d} \citep{kde1d:2024} 
and are represented in blue in the two bottom figures.

The covariance matrix based on the barycentre approach can be estimated by computing the estimates of each matrix in \eqref{Eq:L2_cov_matrixSchur}.
The elements of the main diagonal of $\hat{\bm{\Psi}}$ are the sample estimates of the first moments of the $U_i$. We used the sample means leading to $\hat{\bm{\Psi}}=\diag(-0.66, -0.42, -0.21, -0.21)$. This shows a right-skewed tendency of the latent distributions.
For the computation of $\hat{\bm{\mathfrak E}}_{UU}$, we applied two methods: (i) the elements of the main diagonal, ${\E}(U_i^2)$, were estimated as the sample second moments of the $U_i$; (ii) the elements outside the main diagonal are $\hat{\mathcal E}(U_i ,U_j)=\int_0^1 \hat F_{U_i}^{-1}(t) \hat F_{U_j}^{-1}(t)\ dt$, $i \neq j$, whose integrals were computed using numerical routines from the R package \emph{calculus} \citep{calculus.R:2022}. According to Lemma~\ref{Prp:quantileF}, $\hat F_{U_i}^{-1}(t)=2\hat F_{W_i}^{-1}(t)-1$, where $\hat F_{W_i}^{-1}(t)$ is the estimated quantile function of the fitted distribution, and $W_i \sim {\rm Beta}(\hat a_i, \hat b_i)$, $i=1,2$. For $i=3,4$, $\hat F_{U_i}^{-1}(t)$ was calculated using the function \emph{qkde1d} \citep{kde1d:2024}. This led to the following:
\begin{equation*}
\begin{small}
\hat{\bm{\mathfrak E}}_{UU}=\left[
    \begin{array}{rrrr}
         0.59& \\ 
         0.44& 0.35 \\
         0.35& 0.32& 0.37 \\
         0.34& 0.31& 0.35& 0.34\\
    \end{array}\right].
    \end{small}
\end{equation*}

After estimating the remaining matrices in \eqref{Eq:L2_cov_matrixSchur}, i.e., covariance matrices of the centres and ranges, we obtained the following sample symbolic standard deviations: $10.22$, $15.83$, $75.25$, and $574.45$ for $i=1,\ldots,4$, respectively. Furthermore, the sample correlation matrix is

\begin{equation*}
\widehat{{\rm Cor}}_B{(\bm{X})}=\left[
    \begin{array}{rrrr}
         1.00  \\
         0.85&           1.00   \\
        -0.18&          -0.40&          1.00 \\
        -0.17&          -0.39&          0.99&          1.00\\
    \end{array}\right] .
\end{equation*}

The sample correlation matrix anticipates that departure delays ($x_1$) and arrival delays ($x_2$), as well as time spent in the air ($x_3$) and distance between airports ($x_4$), are highly positively correlated (0.85 and 0.99, respectively). 
The remaining pairs of variables show low to moderate negative correlations. For example, in long-distance flights (or flights of longer duration), it is expected that the pilots can compensate for potential delays, resulting in lower arrival delays. This is expressed by the sample correlation of $-0.39$ ($-0.40$) between $x_2$ and $x_4$ ($x_2$ and $x_3$).

\subsection{RTT dataset}
\cite{salvador2014rtt} introduced a framework to identify traffic redirection attacks, using a group of monitoring probes located across various geographic locations. These probes regularly measured the time it took for a set of 10 data packages to be sent to a target and return, the round-trip time (RTT). The aim of the study was to detect when the data packages relay through a third entity before reaching the target. 
This dataset was fully analysed in \cite{Subtil.et.al:2023}.
The intervals are built from each set of $10$ data packages.  
In this example, we considered the target in Hong Kong, and eight monitoring probes, each corresponding to a variable ($p=8$), located in Amsterdam ($x_1$), Chicago ($x_2$), Vi\~{n}a del Mar ($x_3$), Frankfurt ($x_4$), Hafnarfjordur ($x_5$), S\~{a}o Paulo ($x_6$), and two in Johannesburg, named Johannesburg1 ($x_7$) and Johannesburg2 ($x_8$). We only considered traffic redirected to the Madrid relay, resulting in $n=564$ observations.

At each timestamp, $h$, and probe, $i$, only a few descriptive statistics were recorded, like the sample mean ($\overline{a}_{hi}$) and sample median ($\tilde a_{hi}$) together with the minimum ($a_{h,\min}$) and maximum ($a_{h,\max}$) values of the 10 RTT measures. The limited information about the microdata makes it impossible to fit any distribution to the $U_i$ using the traditional methods. Nevertheless, Pearson’s empirical ``rule of thumb'' allowed us to estimate the mode of each set of microdata, per timestamp and probe: ${\rm mo}_{hi}=3\tilde a_{hi} - 2\bar a_{hi}$, $h=1,\ldots,n$, $i=1,\ldots,p$. The sample means of the modes per probe $\hat m_i=\sum_{h=1}^n mo_{hi}/n$, $i=1,\ldots,p$ are $(-0.14, -0.13, -0.34, -0.58, -0.69, -0.34, -0.17, -0.09)^T$, respectively. This suggests that the distributions of the latent variables are not symmetric. Additionally, the Kruskal-Wallis test concluded that the medians of the probes' modes are different. To test if each latent variable is symmetric we followed the recommendations in \cite[pp.~245]{Sheskin:2011}, and tested if the proportion of modes higher than zero is 1/2 (meaning its median is zero) using the exact test for a binomial proportion (called Clopper-Pearson test), with the Bonferroni correction. The results indicate medians are differ significantly from zero, except for the probes Chicago ($x_2$) and Johannesburg2 ($x_8$). Hence, we considered ${m_2=m_8=0}$. For simplicity, we assumed that the latent variables followed a triangular distribution, i.e., ${U_i \sim {\rm Triang}(-1,1,m_i)}$, where $m_2=m_8=0$, and the other modes, $m_i$, were estimated by $\hat m_i$.

Under the assumption of a triangular distribution, $\E(U_i)=m_i/3$ and $\Var(U_i)=(m_i^2+3)/18$. As a result, $[\hat{\bm{\mathfrak E}}_{UU}]_{ii}=\E(U_i^2)=(m_i^2+3)/6$. Additionally, the quantile function of $U_i$ can be written as
\begin{equation*}
    F_{U_i}^{-1}(t)=\left\{\begin{array}{lc}
         -1+\sqrt{2t(m_i+1)},&  0< t \leq\frac{m_i+1}{2}\\
         1-\sqrt{2(1-t)(1-m_i)},& \frac{m_i+1}{2} < t \leq 1
    \end{array} \right. .
\end{equation*}
Without loss of generality, we can assume $m_i<m_j$, $i\neq j$, and 
\begin{eqnarray*}
   {\mathcal E}(U_i ,U_j)&=&\int_0^{\frac{1}{2}(m_i+1)}  \bigg( -1+\sqrt{2t(m_i+1)}\bigg) \left( -1+\sqrt{2t(m_j+1)}\right) \ dt\\
   &+& \int_{\frac{1}{2}(m_i+1)}^{\frac{1}{2}(m_j+1)}  \bigg( 1-\sqrt{2(1-t)(1-m_i)}\bigg) \left( -1+\sqrt{2t(m_j+1)}\right)\ dt\\
   &+& \int_{\frac{1}{2}(m_j+1)}^1 \bigg( 1-\sqrt{2(1-t)(1-m_i)}\bigg) \left( 1-\sqrt{2(1-t)(1-m_j)}\right) \ dt.
\end{eqnarray*}
For simplicity, the previous integrals were computed numerically using the routines of the R package \emph{calculus} \citep{calculus.R:2022}, leading to the estimated matrix:
\begin{equation*}
\begin{small}
\hat{\bm{\mathfrak E}}_{UU}=\left[
    \begin{array}{cccccccc}
0.170\\ 
0.167& 0.167\\ 
0.168& 0.167& 0.186 \\
0.169& 0.167& 0.170& 0.223\\ 
0.169& 0.168& 0.171& 0.174& 0.247\\ 
0.168& 0.167& 0.169& 0.170& 0.171& 0.186\\ 
0.167& 0.167& 0.168& 0.169& 0.170& 0.168& 0.172\\ 
0.167& 0.167& 0.167& 0.167& 0.168& 0.167& 0.167& 0.167\\
    \end{array}\right] .
    \end{small}
\end{equation*}

Using \eqref{Eq:L2_cov_matrixSchur}, the previously estimated quantities resulted in the sample symbolic standard deviations: $9.44$, $8.80$, $56.35$, $9.97$, $11.50$, $8.13$, $8.37$, and $8.34$, for $i=1,\ldots,8$, and the sample correlation matrix:
\begin{equation*}
\begin{small}
\widehat{{\rm Cor}}_B{(\bm{X})}=\left[
    \begin{array}{rrrrrrrr}
 1.00\\ 
 0.70&  1.00\\ 
 0.15&  0.18&  1.00\\ 
 0.72&  0.85&  0.22&  1.00\\ 
 0.61&  0.72&  0.19&  0.75&  1.00\\ 
-0.02& -0.03& -0.11& -0.08& -0.05&  1.00\\ 
 0.68&  0.78&  0.16&  0.78&  0.66& -0.01&  1.00 \\ 
 0.66&  0.77&  0.16&  0.79&  0.67& -0.05&  0.80&  1.00\\ %
    \end{array}\right] .
    \end{small}
\end{equation*}

Notice that the probe located in São Paulo ($x_6$) exhibits a very small negative linear association with the other probes. This corroborates the findings in \cite{Subtil.et.al:2023}, where the authors suggest that this probe is not useful for the detection of internet attacks. 
Additionally, note that the estimated correlations related to the probe in Viña del Mar ($x_3$) are much smaller in absolute value than the other probes (except São Paulo). The existence of an atypical behaviour of the probe in Viña del Mar was also detected and discussed in \cite{Subtil.et.al:2023}.
The remaining probes, $x_1$, $x_2$, $x_4$, $x_5$, $x_7$, and $x_8$, reveal a moderate ($0.61$) to high ($0.85$) positive correlation among themselves.

To conclude, we remark that even though the results allure to the symmetric triangular distribution, that is, the entries of $\hat{\bm{\mathfrak E}}_{UU}$ look similar to $1/6$ and the main diagonal of ${\hat{\bm{\Psi}} = \diag(0.05,  0, -0.11, -0.19, -0.23, -0.11, -0.06,  0)}$ is close to the zero vector, there is an added value in considering a non-symmetric triangular distribution. The Frobenius norm of the difference between $\widehat{{\rm Cor}}_B{(\bm{X})}$ and the correlation matrix based on $\bm{S}_{CC} + \bm{S}_{RR}/24$ (under the assumption of a symmetric triangular distribution) is $0.174$.
The magnitude of this value accentuates the difference in choosing an assumption that uses partial information about the data as opposed to assuming the symmetric triangular distribution, which is a more simplified approach. 

A common assumption amongst the SDA community is the uniform distribution, which assigns the same level of uncertainty to the microdata. The Frobenius norm of the difference between $\widehat{{\rm Cor}}_B{(\bm{X})}$ and the correlation matrix based on $\bm{S}_{CC} + \bm{S}_{RR}/12$ increases to 0.288.

\section{Conclusions and discussion}\label{Sec:concl}

In this work, we introduced the framework of interval data paired with the distribution of the microdata. Using the model proposed in \cite{Oliveira.et.al:2021}, we started by scaling the microdata in the interval $[-1,1]$, represented by the latent random variable $U_i$. This transformation has the merit of simplifying the theoretical derivations presented. 
 
 Based on it, we deduced explicit formulas for the Mallows' distance between two $p$-dimensional intervals, under very mild assumptions on the distribution of the microdata. In its most general form, the squared Mallows' distance can be decomposed into three terms: the squared Euclidean distance between the two vectors of the centres, the weighted squared Euclidean distance between the two vectors of the ranges, where the weights depend on the second moment of the latent random variables, and a third term that balances the contribution of the centres and ranges, weighted by the expected value of $U_i$. Assuming a symmetric distribution for the latent random variables eliminates this cross-term, turning the squared Mallows' distance into the sum of two squared Euclidean distances: one based on the distance between the centres and the other on the weighted distance between the ranges. In this case, the ranges' weights are quantities in $[0 , 1/4]$, which accentuates the unbalanced contribution of the centres and ranges to the distance.

The general expression also allowed us to argue that the Mallows' distance between two hyperrectangles in $\IR^p$ with the same distribution of the latent random variables in each dimension can be expressed as a weighted Euclidean distance between two points in $\R ^p\times(\R_0^+)^p$ composed by the vector of the centre and range combined. The associated covariance matrix $\bm{H}^{-1}$ is the inverse of a $2$ by $2$ block matrix; each of the blocks is a $p\times p$ diagonal matrix. These block matrices only depend on the first two moments of the latent random variables.

The closed form of the Mallows' distance led to generalising the definitions of the expected value and covariance matrix of an interval-valued random vector. The expected value is defined as the interval that minimises the expected value of the square of the Mallows' distance to the interval-valued random vector called population barycentre or Fréchet population mean. The minimum value of the function to be minimised is the Fréchet variance, which coincides with the trace of the deduced symbolic covariance matrix called the total variance. The deduction of the symbolic covariance matrix, based on the barycentre approach, highlights the contribution of the covariance between centres and ranges $\bm{\Sigma}_{CR}$. This is a novelty, since most of the works in SDA assume a symmetric distribution for the microdata, concealing the role of this matrix.

In practice, we may not have full information about the microdata, and even if we do, it may be difficult to fit a parametric distribution. Our examples illustrate the use of kernel density estimators to overcome this issue. Additionally, we discussed an example where only limited information about the microdata is available.

\backmatter

\bmhead{Acknowledgements}

The authors thank Dr Paulo Salvador and Dr Ana Subtil for sharing the RTT dataset used in one of the examples. This work was supported by Fundação para a Ciência e Tecnologia, Portugal, through the projects [\href{https://doi.org/10.54499/UIDB/04621/2020}{UIDB/04621/2020}, UIDB/04459/2020, UIDP/04459/2020].

\begin{appendices}


\appendix

\section{Proof of results in Section~\ref{Sec:Mallows}}\label{AppendixA}

\subsection{Proof of Lemma~\ref{Prp:quantileF}}\label{AppendixA1}

If $b=0$, then Lemma~\ref{Prp:quantileF} holds trivially, since, for all $t\in(0,1]$, we have ${\rm P}(Z=a)=1$ and $F_Z^{-1}(t)=a$. 
Suppose now that $b>0$. For $z \in \R$,
\begin{equation*}
    F_Z(z) = {\rm P}\Big(a + b~W \leq z\Big) = {\rm P}\Big(W \leq \frac{z-a}{b}\Big) = F_W\Big(\frac{z-a}{b}\Big)\ .
\end{equation*}
By the definition of $F_Z^{-1}(t)$, it follows that
\begin{equation*}
F_Z^{-1}(t) = \text{inf}\{z \in \R:\ t \leq F_Z(z)\}
= \text{inf}\Big\{z:\ t \leq F_W\Big(\frac{z-a}{b}\Big)\Big\} 
= \text{inf}\Big\{a+b~u:\ t \leq F_W(u)\Big\}. 
\end{equation*}
Since $b >0$, we obtain $F_Z^{-1}(t) = a + b ~\text{inf}\{u \in \R:\ t \leq F_W(u)\} = a + b~F_W^{-1}(t),$
as required.

\subsection{Proof of Corollary~\ref{Th:Dist_Wasser2_Geral}}\label{AppendixA2}

By Definition~\ref{def:Mallows_unidim} and Lemma~\ref{Prp:quantileF},
   \begin{align}
        d_M(x_1,x_2)^2\ &=\ \int_{0}^{1}\left(F_{1}^{-1}(t) - F_{2}^{-1}(t)\right)^2 dt\notag \\ 
        &=\ \int_{0}^{1}\left(c_1-c_2 + \frac{ r_1}{ 2}~F_{\U_1}^{-1}(t) - \frac{ r_2}{ 2}~F_{\U_2}^{-1}(t) \right)^2  dt\notag\\
        &=\ (c_1-c_2)^2\ 
        +\ (c_1-c_2) \left( r_1\int_{0}^{1} F_{\U_1}^{-1}(t) \ dt - r_2\int_{0}^{1} F_{\U_2}^{-1}(t)\ dt \right)\notag\\
        &+\ \frac{ 1}{ 4} \int_{0}^{1} \left(r_1~F_{\U_1}^{-1}(t) - r_2~F_{\U_2}^{-1}(t)\right)^2 dt . \label{L2quantile}
\end{align}
Notice that for $k=1,2$,
\begin{equation}
\int_0^1{\Big(F_{\U_j}^{-1}(t)}\Big)^k\ dt\ = 
\int_{\R}{\Big[F_{\U_j}^{-1}\big(F_{\U_j}(u)\big)\Big]^k\ f_{\U_j}(u) \  du}\ = \int_{\R}{u^k\ f_{\U_j}(u)\ du}\ =\ \E{(\U_j^k)}, \label{Eq:E_U2}
\end{equation}
where $f_{\U_j}$ is the probability density function of the absolutely continuous latent random variable $\U_j,\ {j=1,2}$.

\noindent Denoting ${\mathcal E}(\U_1,\U_2) = \int_0^1 {F_{\U_1}^{-1}(t)F_{\U_2}^{-1}(t)}\ dt$, and replacing \eqref{Eq:E_U2} in \eqref{L2quantile}, we obtain  
   \begin{align}
        d_M(x_1,x_2)^2\ &=\ (c_1-c_2)^2\ 
        +\ (c_1-c_2) \left( r_1\E{(\U_1)} - r_2\E{(\U_2)} \right) \notag\\
        &+\ \frac{ r_1^2}{ 4}\E{(\U_1^2)} 
        \ +\ \frac{ r_2^2}{ 4}\E{(\U_2^2)}
        \ -\  \frac{ r_1 r_2}{ 2}~{\mathcal E}{(\U_1 , \U_2)},\label{Eq:dW2.aux1}
\end{align}
concluding the proof.

\subsection{Proof of Corollary~\ref{Cor:Mallows_linked_IrpinoVerde}}\label{AppendixA2.B}

\noindent By adding and subtracting $\left(  r_1\E{(U_1)}/2 -  r_2\E{(U_2)}/2 \right)^2$ to \eqref{Eq:dW2.aux1}, and considering
        $\mu_j = \E({\tilde{V}}_j) = c_j+{r_j}\E(U_j)/2$, and 
        $\sigma_j^2 = \Var({\tilde{V}}_j) = {r_j^2}\Var(U_j)/4$, it follows that
   \begin{align*}
        d_M(x_1,x_2)^2\ &=\ (\mu_1-\mu_2)^2+\sigma_1^2+\sigma_2^2-\frac{ r_1 r_2}{ 2}\left({\mathcal E}{(U_1 , U_2)}-\E{(U_1)}\E{(U_2)}\right)
        \\
        &=\ (\mu_1-\mu_2)^2+(\sigma_1-\sigma_2)^2+2\sigma_1\sigma_2\left(1-\frac{{\mathcal E}{(U_1 , U_2)}-\E{(U_1)}\E{(U_2)}}{\sqrt{\Var{(U_1)}\Var{(U_2)}}}\right)\\
        &=\ (\mu_1-\mu_2)^2+(\sigma_1-\sigma_2)^2+2\sigma_1\sigma_2\left(1-\rho_{12}\right).
\end{align*}

\noindent The last step is shown by using \eqref{transformquantile} in $\rho_{12}$ and simplifying the expressions.

\subsection{Proof of Corollary~\ref{corollary:mahalob}}\label{AppendixA3}

To prove that $\bm{H}^{-1}$ is a covariance matrix, assuming $\Var(U_i)>0$, $i=1,\ldots,p$, we need to show that it is a symmetric positive definite matrix. This is equivalent to showing that the $2p\times 2p$ matrix $\bm{H}$ is itself a symmetric positive definite matrix. Symmetry is easily seen from the definition of $\bm{H}$. It remains to show positive definiteness. 
Let  $\bm{v} = (\bm{v}_1^T,\bm{v}_{2}^T)^T$ be a non-zero real vector, where $\bm{v}_j = (v_{j1},\ldots,v_{jp})^T$, $j=1,2$, and $\bm{v} \neq \bm{0}$. We have
\begin{align}
\bm{v}^T\bm{H}\bm{v}\ &=\ 2\bm{v}_1^T\bm{v}_1\ +\ 2\bm{v}_2^T\bm{\Delta}\bm{v}_2\ +\ 2\bm{v}_1^T\bm{\Psi}\bm{v}_2\ \notag\\
    &=\ 2\sum_{i=1}^p v_{1i}^2\ +\ \frac{1}{2}\sum_{i=1}^{p} {v}_{2i}^2\E(U_i^2)\ +\ 2\sum_{i=1}^p v_{1i}v_{2i}\E(U_i).\label{Eq:HessDP.aux1}
\end{align}
Considering $\E(U_i^2)=\Var(U_i)+\E(U_i)^2$ in \eqref{Eq:HessDP.aux1}, we obtain
\begin{align}
    \bm{v}^T\bm{H}\bm{v}\ 
        &=\ 2\sum_{i=1}^p \left(v_{1i} + \E(U_i)\frac{v_{2i}}{2}\right)^2\ +\ \frac{1}{2}\sum_{i=1}^{p} {v}_{2i}^2\Var(U_i).\ 
    \label{Eq:HessDP.aux3}
\end{align}
We now show that $\bm{v}^T\bm{H}\bm{v}>0$, when $\bm{v}\neq\bm{0}$. The first term in \eqref{Eq:HessDP.aux3} is always non-negative, and the second is strictly positive, for $\bm{v}_2\neq \bm{0}$ (assuming $\Var(U_i)>0$, $i=1,\ldots,p$). Hence, the inequality $\bm{v}^T\bm{H}\bm{v}>0$ holds. 

\noindent Suppose now that $\bm{v}_2 = \bm{0}$, that is, $\bm{v}^T\bm{H}\bm{v} = 2\bm{v}_1^T\bm{v}_1$. 
 Since, by hypothesis, $\bm{v} \neq \bm{0}$, there exists at least one component of $\bm{v}_1$ different from zero, yielding $\bm{v}^T\bm{H}\bm{v}>0$. Hence, $\bm{H}$ is positive definite. 
 
 \noindent Note that without the assumption of positive variance, one can only ascertain that $\bm{H}$ is positive \mbox{semi-definite}.

\section{Proof of results in Section~\ref{Sec4:moments}}\label{AppendixB}

In this section, we prove the results formulated in Section~\ref{Sec4:moments}.

\subsection{Proof of Theorem~\ref{Th:Barycenter.Pop}}\label{Prf:Barycenter.Pop}

We firstly note that the function $g$ defined in \eqref{Eq:Barycenter.Pop} is convex. This is due to the fact that the square of the Mallows' distance is a convex function and the expected value preserves convexity. 
Following \cite{Boyd2004}, the Mallows' distance \eqref{Eq:Mallows_Mahalob} is convex because $\bm{H}$ is a positive semi-definite matrix, as seen in Appendix~\ref{AppendixA3}.

\noindent Now, since $g$ is convex, any critical point is necessarily a global minimum. Hence, we find the points that make the partial derivatives of the objective function relative to $\bm{c}$ and $\bm{r}$ equal to zero. This leads to
\begin{equation}\label{baryproof}
\begin{cases}
\hspace{-2pt}\dfrac{\partial}{\partial\bm{c}} g(\bm{c},\bm{r}) = \bm{0}\\[20pt]
\hspace{-2pt}\dfrac{\partial}{\partial\bm{r}} g(\bm{c},\bm{r})  = \bm{0}
\end{cases} 
\hspace{-15pt}\Leftrightarrow   
\begin{cases}
\hspace{-2pt}-2\E\hspace{-2pt}\left(\hspace{-3pt} \bm{C}-\bm{c} +\dfrac{1}{2}\bm{\Psi} (\bm{R} -\bm{r}) \hspace{-3pt}\right) = \bm{0}\\[15pt]
\hspace{-2pt}-2\E\hspace{-2pt}\left(\hspace{-3pt} \bm{\Delta} ( \bm{R}-\bm{r}) +  \dfrac{1}{2}\bm{\Psi} (\bm{C}-\bm{c}) \hspace{-3pt}\right) = \bm{0}
\end{cases}
 \hspace{-15pt}\Leftrightarrow  
\begin{cases}
 \hspace{-2pt}\E\hspace{-2pt}\left( \bm{C}-\bm{c}\right) = - \dfrac{1}{2}\bm{\Psi} \E\hspace{-2pt}\left(  \bm{R}-\bm{r} \right) \\[15pt]
 \hspace{-2pt}\left(\bm{\Delta} - \dfrac{1}{4}\bm{\Psi}^2 \right) \E\left(  \bm{R} - \bm{r} \right)=\bm{0}
\end{cases}\hspace{-15pt},
\end{equation}
where $\bm{\Delta} - \bm{\Psi}^2/4 =\diag{\left(\Var(U_1),\ldots,\Var(U_p)\right)}/4$. 
According to Definition~\ref{Def:Model_Micro_Macro}, if ${{\rm P}(R_i=0)}=0$, $i=1,\ldots,p$, the latent random variables are absolutely continuous, and the only solution to \eqref{baryproof} is the symbolic hyperrectangle ${\bm{\mu}_B = (\bm{\mu}_C,\bm{\mu}_R, F_{\bm U})}$.

\noindent Consider now the case where some of the components of $\bm{R} = (R_1,\ldots,R_p)^T$ are equal to $0$ with probability $1$.
It follows from Definition~\ref{Def:Model_Micro_Macro} that, for these components, the corresponding latent random variables are degenerate and their variance is zero.

\noindent In case some of the components of $\bm{R} = (R_1,\ldots,R_p)^T$ are equal to $0$ with probability $1$, the respective latent random variables are also degenerate, according to Definition~\ref{Def:Model_Micro_Macro}, and \eqref{baryproof} is trivially satisfied for any range. As such, $\bm{\mu}_B = (\bm{\mu}_C,\bm{\mu}_R, F_{\bm U})$ remains true.

\noindent The minimum value of the objective function \eqref{Eq:Barycenter.Pop}, called Fréchet variance, is a non-negative real number. Thus, we can use the properties of the trace of a matrix to obtain:

\begin{align*}
    V_F(\bm{\mu}_B)\ &=\
    \E\left( (\bm{C}-\bm{\mu}_C)^T(\bm{C}-\bm{\mu}_C) 
        \ +\  (\bm{R}- \bm{\mu}_R)^T\bm{\Delta}(\bm{R}-\bm{\mu}_R) 
        \ +\  (\bm{C}-\bm{\mu}_C)^T \bm{\Psi} (\bm{R}-\bm{\mu}_R)\right) \\
        &=\ \tr\Big(\E\left( (\bm{C}-\bm{\mu}_C)(\bm{C}-\bm{\mu}_C)^T \right)
        \ +\  \bm{\Delta} \E\left( (\bm{R}- \bm{\mu}_R)(\bm{R}-\bm{\mu}_R)^T \right)\\
        &+\   \bm{\Psi}\E\left((\bm{R}-\bm{\mu}_R)  (\bm{C}-\bm{\mu}_C)^T\right) \Big) \\
        &=\ \tr \left(\bm{\Sigma}_{CC} +\bm{\Delta} \bm{\Sigma}_{RR} + \bm{\Psi}\bm{\Sigma}_{RC}\right)\\ 
        &=\ \tr \left(\bm{\Sigma}_{CC} +\bm{\Delta} \bm{\Sigma}_{RR} + \bm{\Sigma}_{CR}\bm{\Psi}\right),
\end{align*}
concluding the proof. 

\subsection{Proof of Corollary~\ref{Crll:PopCov}}\label{Prf:PopCov}

According to Lemma~\ref{Prp:quantileF}, for $i=1,2$, we can write $F_{B_i}^{-1}(t)=\mu_{C_i}+\mu_{R_i}F_{{U}_{Bi}}^{-1}(t)/{2}$, where $\mu_{C_i}=\E(C_i)$ and $\mu_{R_i}=\E(R_i)$. Furthermore, we have  $\E(U_i)=\int_0^1 F_{U_i}^{-1}(t)\ dt $ and ${\mathcal E}(U_1 ,U_2)=\int_0^1 F_{U_1}^{-1}(t) F_{U_2}^{-1}(t)\ dt$. Therefore,
\begin{align*}
    \int_0^1 \left( G_{1}^{-1}(t)-F_{B_1}^{-1}(t)\right) \left( G_{2}^{-1}(t)-F_{B_2}^{-1}(t) \right) dt\ &=\ \left(C_{1}-\mu_{C_1}\right)\left(C_2-\mu_{R_2}\right)\\ 
    &+\ \frac{1}{2}\left(C_{1}-\mu_{C_1}\right)\left(R_2-\mu_{R_2}\right)\E(U_2)\\
    &+\ \frac{1}{2}\left(C_2-\mu_{R_2}\right)\left(R_{1}-\mu_{R_1}\right)\E(U_1)\\ 
     &+\ \frac{1}{4}\left(R_{1}-\mu_{R_1}\right)\left(R_2-\mu_{R_2}\right) {\mathcal E}(U_1 ,U_2),
\end{align*}
whose expected value is the required expression. 
If $U_1$ and $U_2$ are identically distributed, then ${{\mathcal E}(U_1 ,U_2) = \E(U_1^2)}$, and we immediately obtain \eqref{Eq:PopVar}.

\end{appendices}

\bibliography{bibliography}

\end{document}